\documentclass[10pt,a4paper,twoside]{article}
\usepackage{amsthm,amsmath,amssymb}
\usepackage{amsmath}
\usepackage{amsthm}
\usepackage{amssymb}
\usepackage{epsfig}
\usepackage{psfrag}
\usepackage{graphicx}
\usepackage[mathcal]{euscript}
\usepackage{color}

\newtheorem{teor}{Theorem}[section]
\newtheorem{defin}[teor]{Definition}
\newtheorem{lemm}[teor]{Lemma}
\newtheorem{osse}[teor]{Remark}
\newtheorem{prop}[teor]{Proposition}
\newtheorem{defi}[teor]{Definition}
\newtheorem{coro}[teor]{Corollary}
\newtheorem{prob}[teor]{Problem}

\newcommand{\bele}{\begin{lemm}\begin{sl}}
\newcommand{\enle}{\end{sl}\end{lemm}}
\newcommand{\bedef}{\begin{defi}\begin{sl}}
\newcommand{\eddef}{\end{sl}\end{defi}}
\newcommand{\bete}{\begin{teor}\begin{sl}}
\newcommand{\ente}{\end{sl}\end{teor}}
\newcommand{\beos}{\begin{osse}\begin{rm}}
\newcommand{\eddos}{\end{rm}\end{osse}}
\newcommand{\bepr}{\begin{prop}\begin{sl}}
\newcommand{\empr}{\end{sl}\end{prop}}
\newcommand{\bepro}{\begin{prob}\begin{rm}}
\newcommand{\empro}{\end{rm}\end{prob}}
\newcommand{\bede}{\begin{defin}\begin{sl}}
\newcommand{\edde}{\end{sl}\end{defin}}
\newcommand{\beco}{\begin{coro}\begin{sl}}
\newcommand{\enco}{\end{sl}\end{coro}}

\newcommand{\iniziomodifica}{\color{black}}
\newcommand{\finemodifica}{\color{black}}

\newcommand{\iniziomodificag}{\color{black}}
\newcommand{\finemodificag}{\color{black}}

\newcommand{\thspace}{\hspace{3mm}}

\newcommand{\quext}{\quad\text}
\newcommand{\qquext}{\qquad\text}
\newcommand{\de}{\partial}

\newcommand{\Rz}{\mathbb{R}}

\newcommand{\ove}{\overline}
                           
\newcommand{\RR}{\mathbb{R}}

\newcommand{\beeq}[1]{\begin{equation}\label{#1}}
\newcommand{\eddeq}{\end{equation}}

\newcommand{\beeqa}[1]{\begin{eqnarray}\label{#1}}
\newcommand{\eddeqa}{\end{eqnarray}}

\newcommand{\beal}[1]{\begin{align}\label{#1}}
\newcommand{\eddal}{\end{align}}

\newcommand{\bespl}[1]{\begin{split}\label{#1}}
\newcommand{\edspl}{\end{split}}

\newcommand{\bega}[1]{\begin{gather}\label{#1}}
\newcommand{\edga}{\end{gather}}

\newcommand{\beeqax}{\begin{eqnarray*}}
\newcommand{\eddeqax}{\end{eqnarray*}}

\def\qed{\ifmmode 
  \else \leavevmode\unskip\penalty9999 \hbox{}\nobreak\hfill
  \fi
  \quad\hbox{\hskip.5em\vrule width.4em height.6em depth.05em\hskip.1em}}
\def\endproofsym{\qed}
\renewenvironment{proof}[1][Proof]{\trivlist\item[\hskip\labelsep{\hskip0pt
    {\normalfont\scshape#1.}\hskip .321429\parindent}]\ignorespaces}
{\endproofsym\endtrivlist}
\def\endnobox{\def\endproofsym{}\end{proof}\def\endproofsym{\qed}}

\newcommand{\no}{\nonumber}

\newcommand{\beeqao}{\begin{eqnarray}\no}
\newcommand{\bealo}{\begin{align}\no}
\newcommand{\besplo}{\begin{split}\no}
\newcommand{\begao}{\begin{gather}\no}

\newcommand{\nor}[2]{\|#1\|_{#2}}

\newcommand{\duav}[1]{\langle{#1}\rangle}

\newcommand{\+}{\hspace{1pt}}
\newcommand{\perogni}{\forall\,}
\newcommand{\esiste}{\exists\,}
\newcommand{\itt}{\int_0^t}

\newcommand{\io}{\int_\Omega}

\newcommand{\iTT}{\int_0^T}

\newcommand{\iTo}{\iTT\!\io}

\newcommand{\epsi}{\varepsilon}
\newcommand{\ee}{_{\varepsilon}}

\newcommand{\util}{\widetilde{u}}

\newcommand{\bxi}{\boldsymbol{\xi}}

\newcommand{\bbb}{\boldsymbol{b}}

\newcommand{\lhs}{left hand side}
\newcommand{\rhs}{right hand side}

\DeclareMathOperator{\dive}{div}
\DeclareMathOperator{\deriv}{d}

\DeclareMathOperator{\dom}{dom}

\DeclareMathOperator{\Id}{Id}
\DeclareMathOperator{\inte}{int}
\DeclareMathOperator{\sign}{sign}

\DeclareMathOperator{\loc}{loc}

\newcommand{\LDH}{L^2(0,T;H)}

\newcommand{\LIH}{L^\infty(0,T;H)}

\let\TeXchi\chi
\def\chi{{\setbox0 \hbox{\mathsurround0pt
$\TeXchi$}\hbox{\raise\dp0 \copy0 }}}

\newcommand{\alfaciapo}{\widehat{\alpha}}

\newcommand{\betaciapo}{\widehat{\beta}}

\newcommand{\calJ}{{\cal J}}

\newcommand{\calA}{{\cal A}}
\newcommand{\calF}{{\cal F}}

\newcommand{\calB}{{\cal B}}
\newcommand{\calG}{{\cal G}}

\newcommand{\barr}{\overline{r}}

\newcommand{\barO}{\overline{\Omega}}
\newcommand{\barx}{\overline{x}}
\newcommand{\bart}{\overline{t}}

\newcommand{\subr}{\underline{r}}

\newcommand{\isu}{\overline{\iota}}
\newcommand{\igiu}{\underline{\iota}}

\newcommand{\dit}{\deriv\!t}
\newcommand{\dis}{\deriv\!s}
\newcommand{\dix}{\deriv\!x}

\newcommand{\ddt}{\frac{\deriv\!{}}{\dit}}

\newcommand{\ui}{u_\infty}

\newcommand{\ii}{^i}
\newcommand{\iiu}{^{i-1}}
\newcommand{\nnu}{_\nu}

\newcommand{\sumim}{\sum_{i=1}^m}

\newcommand{\pin}{{p_\infty}}
\newcommand{\qin}{{q_\infty}}

\newcommand{\kin}{{\kappa_\infty}}
\newcommand{\lin}{{\ell_\infty}}

\numberwithin{equation}{section}
\begin{document}

\title{\Large Well-posedness and long-time
  behavior for \\ a class of doubly nonlinear equations}
\author{\normalsize Giulio Schimperna\\
{\normalsize  Dipartimento di Matematica, Universit\`a di Pavia,}\\
{\normalsize Via Ferrata 1, I-27100 Pavia, Italy}\\
{\footnotesize e-mail: {\tt giusch04@unipv.it}}\\
\\
{\normalsize Antonio Segatti}\\
{\normalsize Weierstrass Institute for Applied Analysis and Stochastic}\\
{\normalsize Mohrenstrasse 39, 10117 Berlin, Germany}\\
{\footnotesize e-mail: {\tt segatti@wias-berlin.de}}\\
\\
{\normalsize Ulisse Stefanelli}\\
{\normalsize Istituto di Matematica Applicata e Tecnologie Informatiche -- CNR},\\
{\normalsize Via Ferrata 1, I-27100 Pavia, Italy}\\
{\footnotesize e-mail: {\tt ulisse@imati.cnr.it}}
}
\date{}
\maketitle

\begin{abstract}
 This paper addresses a doubly nonlinear parabolic inclusion of the form 
 $$\,\calA(u_t)+\calB(u)\ni f.$$
 Existence 
 of a solution is proved under suitable
 monotonicity, coercivity, and structure
 assumptions on the operators
 $\,\calA\,$ and $\,\calB$, which in particular are both supposed
 to be subdifferentials of functionals on $\,L^2(\Omega)$. 
 Since {\sl unbounded}\/ operators $\,\calA\,$ are included in the analysis,
 this theory partly extends Colli \& Visintin \cite{CV}. 
Moreover, under additional hypotheses
 on $\,\calB$, uniqueness of the solution is proved. Finally,
a characterization
 of $\omega$-limit sets of solutions is given and we 
investigate the convergence of trajectories to  
limit points.
\end{abstract}


  
\vspace{2mm}

\noindent
{\bf AMS (MOS) subject clas\-si\-fi\-ca\-tion:}\thspace
35K55, 35B40.

\vspace{2mm}



\pagestyle{myheadings}
\newcommand\testopari{\sc G.~Schimperna -- A.~Segatti --- U.~Stefanelli}
\newcommand\testodispari{\sc Doubly Nonlinear Equations}
\markboth{\testodispari}{\testopari}


\section{Introduction}
\label{intro2}

The present analysis is concerned with the study of doubly nonlinear parabolic problems of the form
\begin{equation}
  \label{e}
  \alpha(u_t)  - \,\text{div}\, (\bbb (x,\nabla u)) + W'(u) \ni f,
\end{equation}
in a bounded and regular domain $\, \Omega \subset \Rz^d\,$ $\,
(d=1,2,3)$. 
Here $\,\alpha \subset \Rz\times \Rz \,$ and \iniziomodificag
$\, \bbb(x,\cdot):\Rz^d \rightarrow \Rz^d \,$, $x\in \Omega$, 
\finemodificag are assumed to be maximal monotone while $\, W \,$ is a 
\iniziomodifica
$\lambda$-convex function (see the following assumption $({\bf H1})$ and \cite{AGS} for the definition 
of these kind of perturbations of convex functionals) 
\finemodifica
 and $\, W' \,$ is its derivative.

Equations of the type of \eqref{e} stem in connection with phase change phenomena \cite{BDG,BFL,CLSS,FV,LSS,LSS2}, gas flow through porous media \cite{Visintin96}, damage \cite{BS,BSS,FKS,FKNS,MR2}, and, in the specific case $\, \alpha (\lambda r) = \alpha (r) \,$ for all $\, \lambda >0\,$ and $\, r \in \Rz\,$ (which is however not included in the present analysis), elastoplasticity \cite{DalMaso05,Mielke03,Mielke03b,Mielke04}, brittle fractures \cite{DalMaso04}, ferroelectricity \cite{Mielke-Timofte05}, and general rate-independent systems \cite{Francfort-Mielke05,MM,Mi,Mielke-Rossi05, Mielke-et-al02}. 

The applicative interest of \eqref{e} is related to the fact that it may arise in connection with the {\sl gradient flow} of the {\sl energy} functional
\begin{equation}\label{energia}
E(u) = \int_\Omega \big(\phi (x,\nabla u) + W(u) - f u\big) \qquad (D\phi = \bbb),
\end{equation}
with respect to the metric induced by the {\sl dissipation} functional
\begin{equation}\label{dissipazione}
F(u_t) = \int_\Omega \alfaciapo (u_t) \qquad (\partial  \alfaciapo = \alpha),
\end{equation}
where $ \partial \,$ denotes Clarke's gradient \cite{Clarke} in $\, L^2(\Omega)$.
Namely, by taking variations in $\, L^2(\Omega)$, inclusion \eqref{e} may be derived from the kinetic relation
\begin{equation}\label{dinamica}
\partial F(u_t ) + \partial E(u) \ni 0,
\end{equation}
which represents indeed some balance between the variation in energy and the dissipation in the physical system under consideration.

\iniziomodifica
The mathematical interest in \eqref{e} dates back at least to Barbu
\cite{Barbu75} (see also Arai \cite{Arai} and Senba \cite{Senba}) 
and Colli \& Visintin  \cite{CV} obtained \finemodifica several
existence results in an abstract hilbertian framework. 
The reader is also referred to \cite{Co} for some extension to Banach
spaces. In particular, 
the results in \cite{CV} address relation \eqref{e} in the situation
of a linearly bounded function $\, \alpha \,$ 
and a convex potential $\, W$. More recently, Segatti extended some
result of \cite{CV} to classes of $\, \lambda$-convex 
potentials $\, W \,$ \cite{Se}. In the same paper, the existence of a
global attractor is investigated. 
As no uniqueness results are available for \eqref{e} in the genuine
doubly nonlinear case, 
the analysis of \cite{Se} is tailored to Ball's theory of generalized semiflows \cite{Ba1,Ba2}.

The first issue of this paper is indeed to weaken the requirements
of \cite{CV,Se} and prove novel existence and regularity results in
the situation of unbounded graphs $\, \alpha$. 
This extension turns out to be crucial from the point of view of
applications 
since it entails the possibility of choosing dissipation densities $\,
\widehat \alpha \,$ with bounded domain. In fact, a quite common
choice is $\,\widehat \alpha (r) \sim r^2 \,$ for $\, r \geq 0 \,$ and
$\, \widehat \alpha(r)= + \infty \,$ for $\, r <0\,$ representing
indeed the case of irreversible 
evolutions $\, u_t \geq 0$. Let us point out that a first result in
this direction for the case of 
convex potentials $\, W \,$ and linear functions $\, \bbb\,$ has been obtained in \cite{BFL}.

A second main focus of this paper is on uniqueness of solutions to
\eqref{e}. 
Let us stress that uniqueness is generally not expected in genuine
doubly nonlinear situations. 
We however address here the case of locally smooth functions $\, W'\,$
which are singular at 
the bounded extrema of their domain. By requiring some suitable
compatibility between the singular growth of $\, W' \,$ and the
non-degeneracy of $\, \bbb $, we prove that solutions to \eqref{e}
stay uniformly away from the singular points of $\, W'\,$ (separation
property). H
ence, continuous dependence on data follows at least in the specific
case of linear functions $\,\bbb$. 
This is, to our knowledge, the first uniqueness result available for doubly nonlinear equations of the class \eqref{e}.

The third part of the paper is devoted to the description of the
long-time behavior of solutions to \eqref{e}. 
First of all, we prove a characterization of the (non-empty) $\,
\omega$-limit sets of 
trajectories as solutions to a suitable stationary problem. Secondly,
as $\, \bbb = \Id$, $\, W' \,$ is analytic 
in the interior of its domain, and $\, \alpha\,$ fulfills suitable
growth conditions, an application of 
the \L ojasiewicz-Simon \cite{Lo1,Si2} inequality entails that the $\,
\omega$-limit reduces to a 
single point and the whole trajectory converges (this fact turns out
to be particularly interesting 
since the above mentioned stationary limit problem may exhibit a
continuum of solutions). 
Let us stress that both results hold independently of the uniqueness of a solution.

\iniziomodifica
For the sake of completeness, we shall mention that other types of doubly nonlinear equations have attracted a good deal of interest
in recent years. In particular, equations of the form
\begin{equation}\label{altrotipo}
(\calA u)_t+\calB (u)\ni f,
\end{equation}
where $\calA$ and $\calB$ are nonlinear maximal monotone operators
(even nonlocal in time or with an explicit 
time dependence) have been addressed, and various 
existence, uniqueness and long time behaviour results are available. 
With no claim of completeness, we quote 
\cite{Aka,AL,AsFuKe,Ba79,DiBSho,GiSte,Grange-Mignot72,KePa,Mir,Shi,fico,ganzo},
 among many others.
\finemodifica

Plan of the paper. Section \ref{mainres} is devoted to the detail 
of our assumptions and the statement of the main results. Then,
Section~\ref{secexist} 
brings to the proof of the existence result by means of a
time-discretization procedure and passage 
to the limit technique. Section \ref{secseparuniq} is devoted to the
proof of the separation 
property and uniqueness. Finally, Section \ref{seclongtime} addresses the long-time behavior issues.

\section{Statement of the problem and main results}
\label{mainres}

\subsection{Assumptions and preliminary material}
\label{notat}

Let $\Omega$ be a smooth and bounded domain
in $\RR^d$, $d\in\{1,2,3\}$%
\iniziomodifica%
\footnote{We choose $d\in\{1,2,3\}$ just for
simplicity and in view of possible physical applications. However,
 we stress that all the following results hold in any space dimension.}.
\finemodifica
 Set,
for $t>0$, $Q_t:=\Omega\times(0,t)$. Set also
$H:=L^2(\Omega)$ 
and denote by $|\cdot|$ the norm both in $H$ and in $H^d$ and by $\|\cdot\|_X$
the norm in the generic Banach space $X$. Moreover, we
indicate by $(\cdot,\cdot)$ the scalar product in 
$H$ and by $\langle\cdot,\cdot\rangle$ the duality
between a Banach space $V$ and its topological dual $V'$.
In the sequel, the same symbol
$c$ will be used to indicate some positive constants,
possibly different from each other, appearing in the 
various hypotheses and computations and depending only
on data. When we need to 
fix the precise value of one constant, we shall use 
a notation like $c_i$, $i=1,2,\dots$, instead.

In the sequel, we shall present the basic assumptions
and notations which will be kept, with small variations, 
for all the rest of the paper.

\begin{description}
\item{\bf Assumption (H1).} {\sl Let $\beta\subset\RR\times\RR$ be a maximal monotone
graph and define the {\sl domain}\/ of $\beta$ in $\RR$
as the interval 
$I:=\dom_\RR\beta=\{r\in\RR:\beta(r)\not=\emptyset\}$.
We assume, for simplicity, that $\beta$ is normalized 
in such a way that $0\in I$ and 
$\beta(0)=0$. We denote as
$\beta^0$ the {\sl minimal section}\/ 
(cf., e.g., \cite[p.~28]{Br}) of the graph $\beta$.
We shall also indicate by
$\betaciapo:\RR\to[0,+\infty]$
the convex and lower semicontinuous function such that
$\betaciapo(0)=0$ and $\beta=\partial\betaciapo$,
$\de$ denoting the subdifferential in the sense of Convex
Analysis (here in $\RR$). Note that 
\cite[Prop.~2.11, p.~39]{Br}
$D(\betaciapo):=\{r\in\RR:\betaciapo(r)<+\infty\}
\subset\overline I$.
\iniziomodifica
 We assume the potential
$W$
to be $\lambda$-convex,
 namely of the form (since $\widehat{\beta}$ is assumed to be convex)
\beeq{potpos}
  W(r):=\betaciapo(r)-\frac{\lambda\+r^2}2+c_W
   \qquext{for }\,r\in D(\betaciapo),
\end{equation}
where $\lambda\ge 0$ and $c_W$ is an integration 
constant.
\finemodifica
 Moreover, 
we suppose the coercivity property
\beeq{potpos2b}
   \esiste\eta>0:~~W'(r)r\ge\eta r^2
   \qquext{for }\,|r|\,\text{ sufficiently large in }\,\dom_\RR\beta
\end{equation}
(where with some abuse of notation we are
denoting by $W'(r)$ the (multi-)function
$\beta(r)-\lambda r$, see Remark \ref{ultima} below).}
\end{description}
Note that \eqref{potpos2b} entails no restriction if $\, I \,$ is bounded.
Possibly modifying the integration constant $c_W$ in 
\eqref{potpos} it is not restrictive to 
suppose also that, for $\eta$ as above,
\beeq{potpos2c}
   W(r)\ge\frac{\eta r^2}2
   \qquad\perogni r\in D(\betaciapo).
\end{equation}

In the sequel, $\beta$ will be identified with a 
{\sl maximal monotone}\/ operator (still denoted as $\beta$)
from $H$ to $H$, whose {\sl domain}\/ (in $H$) is the set
\beeq{effdom}
  \dom_H\beta:=\big\{v\in H \ \ \text{such that there exists} 
\ \  \eta \in H \ : \ \eta \in \beta(v)\ \  \text{a.e. in} \ \ \Omega\big\}.
\end{equation}

\beos
 \label{puofareinfty}
 We point out that, even if $I$ is open and bounded,
 given $v\in C(\overline\Omega)\cap\dom_H\beta$,
 it is not excluded
 that there exists a nonempty set
 (necessarily of zero Lebesgue measure) 
 $\Omega_\infty$ such that $v(x)\in\partial I$
 when $x\in\Omega_\infty$. 
\eddos
\begin{description}
\item {\bf Assumption (H2).}
{\sl Let $\alpha\subset\RR\times\RR$ be a maximal monotone
graph and suppose that
$0\in J:=\dom_\RR\alpha=\{r\in\RR:\alpha(r)\not=\emptyset\}$.
Moreover, let, for some $\sigma>0$ and $S_-\le0\le S_+$, 
\beeq{alfaprimo}
  \alpha'(r)\ge 2\sigma 
   \qquext{for a.e.~}\/r\in \inte J\setminus [S_-,S_+],
\end{equation}}
\end{description}
Note that $\, \alpha \,$ is single-valued apart (possibly) from
a countable number of points in $\,J\,$ and
that any selection of $\, \alpha \,$ is almost everywhere 
differentiable.
As before, let $\alfaciapo:\RR\to[0+\infty]$
be the convex and lower semicontinuous function such that
$\alfaciapo(0)=0$ and $\alpha=\partial\alfaciapo$.
As above, the operator $\alpha$ will be identified with its
(possibly multivalued) realization on $H$. Clearly,
\eqref{alfaprimo} implies that
\beeq{alfaprimo2}
  \esiste c=c(\sigma,S_-,S_+,\Omega)\ge 0:~~
   \io(\alpha^0)(v)v\,\dix\ge \sigma|v|^2-c 
   \qquad\perogni v\in H.
\end{equation}
\beos\label{ultima}
>From this point on, in order to simplify the presentation,
we will systematically refer to the multi-functions
$\, \beta$, $\, W'$, and $\, \alpha \,$ by using the notation for
(single-valued) functions. 
In particular the symbol $\, \beta(u)\,$
will be used in order to denote some suitable selection 
$\, \eta \in H \,$ such that $\, \eta \in \beta (u)\,$
almost everywhere in $\, \Omega$,
and so on. We believe that this simplification will cause no confusion but 
rather help clarifying some statements and proofs.  This notational
 convention is taken throughout the remainder
of the paper with the exception of distinguished points
where we add some extra comment.
\eddos
\begin{description}
\item {\bf Assumption (H3).} 
{\sl Let $\phi:\Omega\times\RR^d\to[0,+\infty)$ be such
that:
\beal{phi1}
  & \phi(x,\cdot)\in C^1(\RR^d)\quext{for a.e.~}\,x\in \Omega,\\
 \label{phi2}
  & \phi(x,\cdot)\quext{is convex and }\,\phi(x,0)=0~~
  \text{for a.e.~}\,x\in \Omega,\\
 \label{phi3}
  & \phi(\cdot,\bxi)\quext{is measurable~~for all }\,\bxi\in\RR^d.
\end{align}
Then, we can set
\beeq{defib}
  \bbb:=\nabla_{\bxi}\phi: \Omega \times \Rz^d \rightarrow \Rz^d.
\end{equation}
We assume that, for a given $p>1$, $\phi$ satisfies the growth 
conditions
\bealo
  & \esiste \kappa_1,\kappa_2,\kappa_3>0:~~
   \phi(x,\bxi)\ge \kappa_1|\bxi|^p-\kappa_2, \quad
   |\bbb(x,\bxi)|\le \kappa_3(1+|\bxi|^{p-1})\\
  \label{phi4}
  & \mbox{}~~~~~~~~~~~~~~~~~~~~~~~~~~~\text{for a.e.~}\,x\in\Omega
    ~~\text{and all }\,\bxi\in\RR^d.
\end{align}
Finally, we require that at least one of the following properties
holds:
\bealo
  \text{either}\quad &
   p>6/5\text{~~in \eqref{phi4}},\\
 \label{compsub}
   \text{or}\quad &
    \esiste\eta_1>0,~q>2:~~
    W(r)\ge \eta_1r^q-c
   \quext{for all }\,r\in D(\betaciapo).
\end{align}}
\end{description}
Let us now set either $D(\Phi):=H\cap W^{1,p}_0(\Omega)$ 
or $D(\Phi):=H\cap W^{1,p}(\Omega)$ depending whether homogeneous
Dirichlet or homogeneous
Neumann boundary conditions are taken into account.
Both spaces $W^{1,p}_0(\Omega)$ and  $W^{1,p}(\Omega)$
are intended to be endowed with the full
$W^{1,p}$-norm in the sequel.
Then, we introduce the functional 
\beeq{Phi}
  \Phi:H\to[0,+\infty],\quad
   \Phi(u):=\io\phi(x,\nabla u(x))\,\dix,
\end{equation}
which is defined as identically
$+\infty$ outside its domain $D(\Phi)$.

By means of \cite[Thm.~2.5, p. 22]{Gi}
one has that $\Phi$ is convex and lower semicontinuous on $H$.
Thus, by \cite[\S~2.7]{Br}, its subdifferential 
$\partial\Phi$ in $H$ is a maximal monotone operator in $H\times H$. 
Moreover, we can introduce the operator
$B:\dom_H B\to H$ by setting
\beeq{Bscal}
  (Bv,z):=\io \bbb(x,\nabla v(x))\cdot\nabla z(x)\,\dix,
\end{equation}
for all $z\in D(\Phi)$, defining
\beeq{domB}
  \dom_H B:=\left\{v\in D(\Phi):  
\sup_{ z \in D(\Phi)\setminus\{0\}} { \left|\io 
\bbb(x,\nabla v(x))\cdot\nabla z(x)\,\dix\right|}/{|z|} < + \infty\right\},
\end{equation}
and extending by continuity without changing symbols to $z\in H$. 
In the case of Dirichlet conditions, the latter positions simply 
mean that
$Bv=-\dive(\bbb(\cdot,\nabla v(\cdot))$, where the 
divergence is intended in the sense of distributions.

The inclusion $B\subset\partial\Phi$ is almost immediate to check. 
It is however remarkable that also the converse inclusion can be proved 
(see \cite[Ex.~2.4]{NSV}) by showing that
$B$ is indeed maximal. We shall not provide here a direct proof
of this fact but rather sketch the main steps toward its check.
First of all, one could consider some suitable approximation   
${\boldsymbol b}_\epsi$, $\epsi >0$,
(by a careful mollification, for instance) of the function ${\boldsymbol b}$
and solve the maximality problem for the regularized operator.
Then, by estimating a priori and extracting suitable sequences of indices, 
we prove
that some sequences of solutions to the regularized problem
admit a limit as $\epsi \rightarrow 0$.
Finally, by exploiting the maximality of 
$ {\boldsymbol \xi}(\cdot) \mapsto {\boldsymbol b}(\cdot, {\boldsymbol \xi}(\cdot))$ 
from $L^p(\Rz^d)$ to $L^{p'}(\Rz^d)$, $\, p' = p/(p-1)\,$ (we used here \eqref{phi4}), we identify the limit and
conclude.

Defining now the energy $E:H\to [0,+\infty]$ as
\begin{equation}
  \label{someenergy}
  E(u):= \int_\Omega \big(\phi(x,\nabla u(x)) + W(u(x))\big)\,\dix,
\end{equation}
it is not difficult to see
that $E$ is a \iniziomodificag\
quadratic \finemodificag\ perturbation 
(due to the $\lambda$-term in \eqref{potpos})
of a convex functional. Moreover, we
have
\bele\label{lemmacompsub}
 $E$ has {\sl compact sublevels}\/ in $H$.
\enle
\begin{proof}
 Let $\{v_n\}\subset H$ a sequence of bounded energy.
 If $p>6/5$, $\{v_n\}$ is precompact in $H$ thanks to
 \eqref{phi4} and compactness of the embedding
 $W^{1,p}(\Omega)\subset H$. If, instead,
 the second of \eqref{compsub} holds, a priori
 $\{v_n\}$ only admits a subsequence which is strongly converging
  in $L^p(\Omega)$ and bounded in $L^q(\Omega)$. 
However, since this implies,
 up to a further extraction,
 convergence almost everywhere, the conclusion follows.
\end{proof}
The next simple 
Lemma describes the structure of
$\partial E$.
%
\bele\label{lemmamaxmon}
 The subdifferential of $E$ coincides with the sum
 $\calB:=B+\beta-\lambda\Id$, with 
 domain $\dom_H\calB=\dom_H B\cap\dom_H\beta$.
\enle
\begin{proof}
 It is enough to work on the nonlinear
 part of $\calB$, i.e.~to show that 
 \beeq{soloconv}
   B+\beta=\partial\Big(E+\lambda\frac{|\cdot|^2}2\Big).
 \end{equation}
 The inclusion $\subset$ is clear. To get the converse one,
 it suffices to show that $B+\beta$ is maximal. With this
 aim, we proceed by regularization and replace
 $\beta$ with its (Lipschitz continuous) Yosida approximation 
 \cite[p.~28]{Br} $\beta\ee$, with $\epsi>0$ intended to
 go to $0$ in the limit.
 Then, given $g\in H$,
 by \cite[Lemme 2.4, p.~34]{Br}
 the elliptic problem 
 \beeq{probepsi}
   v\ee\in H,\qquad
    Bv\ee+\beta\ee(v\ee)+v\ee= g,
 \end{equation}
 has at least one solution $v\ee\in H$ which belongs, 
 additionally, to $\dom_H B$. In other words, the operator
 \beeq{Bbetaee}
   B+\beta\ee:\dom_H B\to H
 \end{equation}
is maximal monotone. As $\beta\ee$ is Lipschitz continuous and 
$\beta_\epsi(0)=0$,
it is possible to test \eqref{probepsi}
by $\beta\ee(v\ee)$, obtaining,
for some $c>0$ independent of $\epsi$,
\beeq{stee1}
  |\beta\ee(v\ee)|+|Bv\ee|\le c.
\end{equation}
The boundedness of $\beta\ee(v\ee)$ in $H$
permits then to apply \cite[Thm.~2.4, p.~34]{Br}
which concludes the proof.
\end{proof}
\begin{description}
\item {\bf Assumption (H4).}
{\sl Let us be given
\bega{hypf}
\iniziomodifica
f \in W^{1,1}_{\loc}([0,+\infty);H),\\
\finemodifica
 \label{hpu0}
  u_0\in \dom_H\calB.
\end{gather}
}
\end{description}


\subsection{Main Theorems}

We shall state our abstract Cauchy problem as follows
\bega{eqn}
  \alpha(u_t)+B u+W'(u)= f, \qquext{in }\,H,\\
 \label{iniz} 
  u_{|t=0}=u_0, \qquext{in }\,H.
\end{gather}
Our first result is concerned with the existence of at least one solution to the 
above problem. Its proof will be outlined in Section~\ref{secexist}.
\bete\label{teoexist}{\bf [Existence]}
 Let us assume\/ {\rm ({\bf H1})--({\bf H4})} 
 and take $T>0$. Then, 
 there exists a function $u:Q_T\to\RR$ such that
 \beal{regu1}
   & u\in W^{1,\infty}(0,T;H), \quad
   Bu,~\beta(u)\in L^\infty(0,T;H),\\
  \label{regu2}
   & E(u)\in L^\infty(0,T),
 \end{align}
 which satisfies\/ \eqref{eqn} a.e.~in~$(0,T)$ and the Cauchy
 condition~\eqref{iniz}. More precisely,
 there exists $c>0$ depending on 
 $u_0$, $f$, and $T$ such that
 \beeq{reguunif}
   \nor u{W^{1,\infty}(0,T;H)}
    +\nor{Bu}{L^\infty(0,T;H)}
    +\nor{\beta(u)}{L^\infty(0,T;H)}
    +\nor{E(u)}{L^\infty(0,T)}\le c.
 \end{equation}
 Finally, if either 
 \bega{newrego1a}\tag{{\bf S0}}
   S_-=S_+=0\text{~~in~\,\eqref{alfaprimo}},\\
  \label{newrego1b}\tag{{\bf f1}}
   f\in L^2(0,+\infty;H),~~~f_t\in L^1(0,+\infty;H)
 \end{gather}
 or
 \beeq{newrego2}\tag{{\bf f2}}
   f\in H
    \quext{is independent of time},
 \end{equation}
 then the constant $c$ in\/ \eqref{reguunif} is independent of $T$.
\ente
\beos
Referring to the previous Remark \ref{ultima}, we note that, in case $\, \alpha\,$ and $\, \beta \,$ are multi-valued, the statement above has to be intended in the following sense: we 
prove the existence of selections $\, \xi, \, \eta \in L^\infty(0,T;H)\,$ 
such that $\, \xi \in \alpha(u_t) \,$ and $\, \eta \in\beta(u) \,$
almost everywhere in $\, \Omega \times (0,T) \,$ and 
\begin{equation}
  \label{eq:LO}
  \xi + Bu + \eta - \lambda u =f \qquad \text{in} \ \ H, \ \ \text{a.e. in} \ \ (0,T).
\end{equation}
\eddos

Let us make precise
the comparison of our result with the 
corresponding Theorem proved in \cite{CV}.  Indeed, in our notation, 
\cite[Thm.~2.1]{CV} 
(which is, among the various results of 
\cite{CV}, the closest to our theory)
can be formulated as follows:
\bete{\bf [Colli \& Visintin, 1990]}\label{teoCV}
 Let\/ {\rm ({\bf H1})--({\bf H3})} hold with $\lambda=0$
 in\/ \eqref{potpos}, let $T>0$,
 and let  $\alpha$ be\/ {\rm sublinear}, namely:
 \beeq{alfasublin} 
   |\alpha(r)|\le c(1+|r|)\qquext{for all }\,r\in\RR
    \quext{and some }\,c>0.
 \end{equation}
 Moreover, in place of\/ {\rm ({\bf H4})}, assume
 $E(u_0)<\infty$
 and $\, f \in L^2(0,T;H)$.
 Then, \eqref{eqn}--\eqref{iniz} admits at least one solution
 $u$ satisfying\/ \eqref{regu2} together with
 \beeq{reguCV}
   u\in H^{1}(0,T;H), \qquad
    \alpha(u_t),~Bu,~\beta(u)\in \LDH.
 \end{equation}
\ente
By comparing Theorems \ref{teoexist} and \ref{teoCV} we readily 
check that our result refers to more general classes of nonlinearities
$\, \alpha \,$ and provides more regular solutions. On the other
hand, in place of our assumption \eqref{hpu0}, in \cite{CV} the authors
were able to consider less regular initial
data, satisfying the {\sl finite energy} condition
$E(u_0)<\infty$.

One of the main novelties of the present work is that for a
significant class of potentials $W$ we are able to show
uniqueness and further regularity. 
In order this to hold we have to introduce the following.

\begin{description}
\item {\bf Assumption (H5).}
{\sl Suppose that, for some $\nu>0$, one has
 \beeq{eqsep} 
   \dom_H B \subset C^{0,\nu}(\overline\Omega)
 \end{equation}
 and that there exists a nondecreasing
 function $\gamma:[0,+\infty)\to[0,+\infty)$ such that
 \beeq{eqsepbis}
   \nor{v}{C^{0,\nu}(\barO)}
    \le \gamma\big(|v|+|Bv|\big)
    \qquad\perogni v\in\dom_H B.
 \end{equation}
 Moreover, let $I=\dom_\RR\beta$ be an 
 {\sl open}\/ set and, in case $I$ is, e.g.,
 right-bounded, set $\barr:=\sup I$
 and assume in addition
 that there exist $c>0$, $r_1<\barr$, 
 and $\kappa>0$ such that 
 \beeq{crescibeta} 
   \beta^0(r)\ge \frac c{(\barr-r)^\kappa},
    \qquad\perogni r\in(r_1,\barr), 
 \end{equation}
 where $\kappa$ and $\nu$ satisfy the compatibility condition
 \beeq{kappa}
   2 \kappa\+\nu\ge d.
 \end{equation}
 Analogously, if $I$ is left-bounded,
 set $\subr:=\inf I$ and assume 
 the analogous of \eqref{crescibeta} in a right neighborhood
 $(\subr,r_0)$ of $\subr$, together with \eqref{kappa}.}
\end{description}

\beos\label{cosesepa}
 Condition \eqref{crescibeta} states that $\beta$ explodes
 at $\barr$ at least as fast as a sufficiently high
 negative power of $(\barr-r)$, made precise by 
 \eqref{kappa}. 
 Of course, the higher is this power, the faster the system
 is forced by energy dissipation to keep the solution
 away from the barrier at $r=\barr$. Note that, unfortunately,
 even for $d=1$, this
 excludes from our analysis the {\sl logarithmic}\/
 potentials of the type
 \beeq{logpot}
   W(r)=(1-r)\log(1-r)+(1+r)\log(1+r),
 \end{equation}
 which are relevant for applications
 (cf, e.g.,~\cite{MZ}), since
 they are not enough coercive in proximity of the
 barriers (here given by 
 $\subr=-1$, $\barr=1$). However 
 (cf.~the forthcoming Remark \ref{quandosepa}), in one dimension 
($d=1$) it 
 is possible to consider some potentials $W$ 
 which are {\sl bounded}\/ on $[\subr,\barr]$
 and thus qualitatively similar to \eqref{logpot}.
\eddos
\beos\label{quandosepa}
 Property \eqref{eqsep} depends of course on the operator $B$
 and on the space dimension $d$. For instance, if $d=3$ and 
 $B=-\Delta$, then \eqref{eqsep} holds for $\nu\leq 1/2$.
 If $B$ is a nonlinear operator coming from a functional 
 $\Phi$ of $p$-growth (cf.~\eqref{phi4}), $p\ge2$,
 satisfying the regularity conditions 
 \cite[(22), (35)]{Sa} (this is the case, e.g., of the $p$-Laplace operator), and moreover
 Dirichlet conditions are taken 
 (i.e.~$D(\Phi)$ is chosen as $W^{1,p}_0(\Omega)$), 
 then one can apply, e.g., \cite[Thm.~2]{Sa}, yielding
 $\dom_H B\subset W^{\zeta,p}(\Omega)$ for
 any $\zeta<1+1/p$, which in turn
 entails \eqref{eqsep}, \iniziomodificag
 with a suitable $\nu$, \finemodificag
 for $p\ge2$ if $d=2$ and for $p>2$
 if $d=3$.
\eddos
\bepr{\bf [Separation property]}\label{teosep}
 Let us assume\/ {\rm ({\bf H1})--({\bf H5})}
 and let $u$ be any solution in the regularity setting
 of\/ {\rm Theorem~\ref{teoexist}}.
 Then, $u(t)\in C^{0,\nu}(\overline\Omega)$ for\/
 {\rm every} $t\in [0,T]$ and 
 $\esiste c_0>0$ such that
 \beeq{reguinfty}
   \nor{u(t)}{C^{0,\nu}(\overline\Omega)}
    \le c_0 \qquad\perogni t\in[0,T].
 \end{equation}
 Moreover, $\esiste \subr_*<0,\barr^*>0$, with 
 $\barr^*<\barr$ if $I$ is right-bounded and
 $\subr_*>\subr$ if $I$ is left-bounded, 
 such that
 \beeq{sepsu}
    \subr_*\le u(x,t)\le \barr^*\qquad\perogni (x,t)\in \overline Q_T.
 \end{equation}
 Finally, $c_0$, $\subr_*$, and $\barr^*$ are 
 independent of $T$ if either
 \eqref{newrego1a}--\eqref{newrego1b}
 or\/ \eqref{newrego2} hold.
\empr
%
%
Let us now come to uniqueness, which requires the linearity of $\, B $, 
i.e.~the following
additional assumption.

\begin{description}
\item {\bf Assumption (H6).}
{\sl Let us consider a matrix field $D$ such that 
\beal{defiD}
 & D\in L^\infty(\Omega;\RR^{d\times d}),
   \qquad D(x)\ \ \text{is symmetric for a.e.}\ \ 
x\in\Omega,\\
 \label{Dpos}
 & D(x)\bxi\cdot\bxi\ge a|\bxi|^2
   \qquext{for all }\,\bxi\in\RR^d\ \ \text{and a.e.} \ \ x\in\Omega \ \ \text{and some} \ \ a>0.
\end{align}
Then, we assume that
$\phi(x,\bxi):=\frac12(D(x)\bxi)\cdot\bxi$
(cf.~({\bf H2})) and define $\Phi$ and $B$
correspondingly, with the choice of 
$D(\Phi)=H^1(\Omega)$, or $D(\Phi)=H^1_0(\Omega)$, 
depending on the desired boundary conditions.}
\end{description}
Setting $V:=D(\Phi)$, which is now a Hilbert
space, $B$ turns out to be linear elliptic
from $V$ to $V'$ and satisfies, for $a$ as in
\eqref{Dpos},
\beeq{Dpos2}
  \duav{Bv,v}\ge a|\nabla v|^2, \quad
   \duav{Bv,z}\le \|D\|_{L^\infty(\Omega;\Rz^{d \times d})}|\nabla v|\+|\nabla z|,
    \qquad\perogni v,z\in V.\\[2mm]
\end{equation}
\bete{\bf [Uniqueness]}\label{uniq}
 Let us assume\/ {\rm ({\bf H1})--({\bf H6})}
 and\/ \eqref{newrego1a}. Moreover, 
 suppose that 
 \beeq{betalip}
   \beta\/\text{ coincides with a locally Lipschitz continuous 
    function in }\  I.
 \end{equation}
 Then, the solution $u$ provided by\/
 {\rm Theorem~\ref{teoexist}}
 is unique.
\ente
\beos\label{dopouniq}
Let us observe that 
 \eqref{betalip} excludes the presence of vertical segments
 in the graph $\beta$, but does not give any further 
 restriction on its behavior at the boundary of $I$.
 As regards conditions \eqref{eqsep}, \eqref{eqsepbis}
 in ({\bf H5}), we note that, whenever ({\bf H6}) holds, these 
 are related to the regularity properties of $D$
 and the space dimension. For instance, if 
 $D$ is Lipschitz continuous, then it is well known 
\cite[Thm. 8.12, p. 186]{Gilbarg-Trudinger} 
that $\dom_H B\subset H^2(\Omega)$
 (recall that we are always supposing $\Omega$ sufficiently
 smooth). Instead, if it is only
 $D\in L^\infty(\Omega)$, then  
 $\dom_H B\subset H^\zeta(\Omega)$
 for any $\zeta<3/2$ \cite[Rem. 4.4]{Sa}.
\eddos
The next results are related to 
existence and characterization of $\omega$-limit
sets of {\sl trajectories}\/ of the system. 
We consider here global solutions $u:[0,+\infty)\to H$
to \eqref{eqn}--\eqref{iniz} satisfying the
regularity frame of Theorem~\ref{teoexist}.
We remark once again that what follows 
holds for any such function $u$ and 
does not require uniqueness.
\bete{\bf [Characterization of the $\,{\boldsymbol \omega}$-limit set]}\label{longtime}
 Let\/ {\rm ({\bf H1})--({\bf H4})}, \eqref{newrego1a}, and 
 either\/ \eqref{newrego1b}
 or \eqref{newrego2} hold.
 Then, for all functions $u$ whose existence is stated in 
{\rm Theorem \ref{teoexist}} and any sequence of times $\{t_n\}$
 with $t_n\nearrow+\infty$, 
 there exist a not relabeled subsequence
 and a function $\ui\in H$  such that  
 \beeq{coinfty1}
    u(t_n) \to \ui 
   \quext{strongly in }\,H.
 \end{equation}
 Moreover, $\ui$ is a solution of the stationary problem
 \beeq{prostaz}
   B\ui+W'(\ui)= g \quext{in } \,H,
 \end{equation}
 where $g=0$ if\/ \eqref{newrego1b}
 holds and $g=f$ if, instead, \eqref{newrego2} is satisfied.
\ente
\beos\label{quandonon0}
In the \iniziomodifica theorem \finemodifica above we are implicitly assuming that the multifunction $\, \alpha \,$ is such that 
\beeq{alpha00}
  \alpha(0)=\{0\}.
  \end{equation}
 It will be clear from the proof that, if \eqref{alpha00}
 does not hold, then we still have \eqref{coinfty1}, 
 but in place of \eqref{prostaz} we have the weaker
 $-B\ui-W'(\ui)+g\in \alpha(0)$.
Moreover, let us point out that our analysis could be easily extended in order to cover the case $\, f = f_1 + f_2 \,$ where $\, f_1\,$ and $\, f_2 \,$ fulfill \eqref{newrego1b} and \eqref{newrego2}, respectively.
\eddos
It is well known that, since $W$ is not convex, problem \eqref{prostaz}
may well admit infinite solutions \cite{Haraux}. 
Thus, the question of convergence
of {\sl all}\/ the trajectory $u(t)$ to one of these solutions is 
a nontrivial one and is not answered by the preceding Theorem. Here,
we are able to show this property under more restrictive 
assumptions, and the basic tool in the proof is the 
so-called \L ojasiewicz-Simon method.
To detail this technique, originally devised in \cite{Lo1,Lo2,Si2},
let us state and comment the further assumptions we need.

\begin{description}
\item {\bf Assumption (H7).}
{\sl Assume that, if \eqref{newrego2} holds, then in
addition 
\beeq{ultregof}
  f\in L^\infty(\Omega), \qquad
   |f(x)|\le M\quext{for a.e.~}\/x\in\Omega
\end{equation} 
and some $M>0$. Moreover, assume that
\bealo
  & \text{there exists an open interval }\,I_0
   =(\igiu,\isu),~~\text{with }\,
   [\igiu,\isu]\subset I,\\
  \label{Iquasiaperto}
  & \text{such that }\,W'(r)+M<0~~\perogni r\le\igiu,\quad
   W'(r)-M>0~~\perogni r\ge\isu,
\end{align}
where $r$ is taken in $I$. In case \eqref{newrego1b}
holds, assume the same with $M=0$. Suppose also that
\beeq{analit}
   W_{|I_0} \quext{is {real analytic}}.
\end{equation}}
\end{description}
\beos\label{qualisolu}
 If ({\bf H7}) holds, a simple maximum principle argument shows that 
 there exists $\epsilon>0$ such that 
 {\sl any}\/ solution $\ui$ to \eqref{prostaz} fulfills
 \beeq{uisepar}
    \igiu+\epsilon\le\ui(x)\le\isu-\epsilon \quad\perogni x\in\Omega.
 \end{equation}
 Note also that it is not excluded that $[\igiu,\isu]=I$.
 This means that $\beta$, and hence $W'$,  
 may contain vertical half lines 
 at the extrema. However, $W'-g$ must
 {\sl have the right sign} at least in a neighborhood
 of the barriers. A similar class of potentials 
 has been considered in \cite[Thm.~2.7]{GPS1}.
\eddos
Let us now choose $B=-\Delta\,$ along with either homogeneous Dirichlet
or Neumann boundary conditions, and take $V$ 
correspondingly (see
({\bf H6})). Consider
again a solution $\ui$ to \eqref{prostaz}. 
Then, the \L ojasiewicz-Simon inequality 
(cf.~\cite[Prop.~6.1]{FS}, see also~\cite{AF}) 
states that 
there exist $c_{\ell},\epsilon>0$, $\theta\in(0,1/2)$
such that
\beeq{Lojfeireisl}
  |E(v)-E(\ui)|^{1-\theta}
   \le c_\ell\|Bv+W'(v)\|_{V^*}
\end{equation}
for all $v\in V$ such that 
\beeq{cosaserveperloj}
   \|v-\ui\|_V\le\epsilon.
\end{equation}
\beos\label{dirneum}
 Inequality \eqref{Lojfeireisl} has been shown in
\cite[Thm. 2.1]{Haraux-Jendoubi99}
 in the Dirichlet case. The analogous statement for
 the Neumann case can be found in
 \cite[Prop. 2.4]{AF}.
\eddos
\beos\label{altriB}
 It should be possible to prove the same result for more general
 (linear and symmetric) operators $B$ satisfying ({\bf H6}). Namely,
 in order to extend 
 inequality \eqref{Lojfeireisl} (or the corresponding {\sl local}\/
 version proved in \cite{FS,AFIR}), what seems to be needed
 is that the coefficients of $B$ are so regular that
 any solution $v$ to the equation $v+Bv=h$ for 
 $h\in L^\infty(\Omega)$ lies
 in $W^{2,p}$ for $p>d$ (cf.~also \cite[Sec.~3]{Je}).
\eddos
We can now state our characterization result for the
$\omega$-limit.
\bete{\bf [Convergence to the stationary state]}\label{Loj}
 Assume\/ {\rm ({\bf H1})--({\bf H4})}, 
 \eqref{newrego1a}, either\/
 \eqref{newrego1b} or \eqref{newrego2},
 and\/ {\rm ({\bf H7})}. Take also 
 $B=-\Delta$.
 Moreover, if \eqref{newrego1b} holds, assume also
 that there exist $c,\xi>0$ such that 
 \beeq{decadef}
   t^{1+\xi}\int_t^\infty |f(s)|^2\,ds \le c 
    \quext{for all }\,t\ge0.
 \end{equation}
 Finally, assume $J=\dom_\RR\alpha=\RR$ and that
 there exist $\sigma',\kappa_\infty,\ell_\infty>0$ and\/
 $1\le p_\infty\le q_\infty$ 
 such that, for all $r\in\RR$,
 \beeq{alfapotenza}
   \sigma|r|+\kappa_\infty|r|^{p_\infty}
    \le|\alpha(r)|
    \le\sigma'|r|+\ell_\infty|r|^{q_\infty}.
 \end{equation}
 Finally, assume the constraint
 \begao
   \chi q_\infty \le (p_\infty+1),
    \quext{where}\\
  \label{relapq}
   \chi\in[1,2]~~\text{if }\,d=1,\quad
   \chi\in(1,2]~~\text{if }\,d=2,
    \quad\text{and }\,\chi\in [6/5,2]~~\text{if }\,d=3.
 \end{gather}
 Then, letting $u$ be a solution, 
 the $\omega$-limit of\/ $u$
 consists of a\/ {\rm unique}
 function $\ui$ solving \eqref{prostaz}, where $g$ 
 is as in\/ {\rm Theorem~\ref{longtime}}.
 Furthermore, as $t\nearrow+\infty$,
 \beeq{coinfty2}
    u(t) \to \ui \qquext{strongly in }\,V\cap C(\barO),
 \end{equation}
 i.e., we have convergence 
 for the\/ {\rm whole trajectory $u(t)$}.
\ente
\beos\label{rempreco}
 Let us point out that, as one proves
 convergence of the trajectory with
 respect to the norm of $H$, then \eqref{coinfty2}
 is immediate since, by \eqref{reguunif}
 and being $B=-\Delta$, the trajectory 
 associated to $u$ is precompact in $V\cap C(\barO)$.
 The same might be extended 
 to more general $B$ satisfying
 ({\bf H6}) (cf.~Remark~\ref{altriB}).
\eddos
\beos\label{qualipq}
 Note that the assumptions entail $\qin\le 5$ if 
 $d=3$. In general, the meaning of \eqref{relapq} 
 is that functions $\alpha$
 with linear growth at $0$ and with behavior at $\infty$ 
 ranging between two not too different powers 
 $\pin\le\qin$ are admissible.
\eddos
\beos\label{suloj}
 The \L ojasiewicz-Simon method used here
 is a powerful tool permitting to 
 prove convergence of trajectories associated 
 to various types of either parabolic or damped 
 hyperbolic evolution equations (we quote, among
 the many recent works, \cite{AF,AFIR,CJ,FS,Je}). 
 In particular, relation
 \eqref{Lojfeireisl} is a {\sl local}\/
 version of the \L ojasiewicz-Simon inequality,
 which has been proved 
 to be useful to treat energy functionals
 of the form \eqref{someenergy}.
 Many other versions of the inequality are 
 available in the recent literature.
\eddos
\beos\label{altrioper}
 In the Theorem above, we just considered
 operators $B$ satisfying {\rm ({\bf H6})}. 
 The main reason for this restriction is that,
 up to our knowledge, inequalities
 of  \L ojasiewicz-Simon type 
 like \eqref{Lojfeireisl} are known
 only for {\sl linear and symmetric}\/ elliptic
 operators. If one could prove \eqref{Lojfeireisl}, or a 
 generalization of its, for other classes of
 elliptic operators, then Theorem~\ref{Loj}
 may be extended as well. 
\iniziomodifica
In this regard, we refer to
 \cite{CF} where a first attempt to include nonlinear elliptic 
operators in the framework of the \L ojasiewicz-Simon inequality is made.
\finemodifica
\eddos


\section{Proof of existence}
\label{secexist}

Let us show here Theorem~\ref{teoexist} by means of a time-discretization and passage to the limit procedure.
\subsection{Regularization and time discretization}
Here we partly follow \cite[Sec.~3]{LSS2}. 
First of all, we substitute the graphs
$\alpha$ and $\beta$ in \eqref{eqn} by their Yosida regularizations
$\alpha\nnu$ and $\beta\ee$, with $\nu,\epsi>0$ both intended to
go to $0$ in the limit. Correspondingly, we set 
$W\ee'(r):=\beta\ee(r)-\lambda\+r$ for $r\in\RR$.

\iniziomodifica
Next, we implement a semi-implicit discretization by means of the backward Euler scheme. To this end, 
we let $N$ be a positive integer and uniformly subdivide $[0,T]$ into subintervals of lenght $\tau:=T/N$.
\finemodifica
We also put $u^0:=u_0$ and discretize $f$ as follows:
\beeq{fi}
  f^i:=f(i\tau)\in H
   \qquext{for }\,i=1,\dots,N,
\end{equation}
the latter position being admissible since $\, f \,$ is continuous in $\, H\,$ (cf.~\eqref{hypf}).
Next, still for $i=1,\dots,N$,
we consider the following scheme
\beeq{eqi}
  \alpha\nnu\Big(\frac{u\ii-u\iiu}{\tau}\Big)
   +B u\ii+\beta\ee(u\ii)+\epsi u\ii
   =\lambda u\iiu+f\ii.
\end{equation}
By induction on $i$, it is easy to show existence 
for this scheme. Actually, as $i$ is fixed, the \rhs\ of \eqref{eqi} is 
a known function in $H$. Moreover, the operator
\beeq{sommagra}
  {\cal D}_{\epsi,\nu}:H\to H, \qquad
   {\cal D}_{\epsi,\nu}:v\mapsto\alpha\nnu\Big(\frac{v-u\iiu}{\tau}\Big)
   +\beta\ee(v),
\end{equation}
is maximal monotone and Lipschitz continuous. In particular,
its domain is the whole space $H$. Since also $B$ is maximal 
monotone, we can apply, e.g., \cite[Cor.~2.7, p.~36]{Br},
which yields that the sum $B+{\cal D}_{\epsi,\nu}$ is also maximal
monotone. Hence, the operator 
$\epsi\Id+B+{\cal D}_{\epsi,\nu}$ acting on $u^i$ on the 
\lhs\ of \eqref{eqi} is onto, which 
means that we have existence of a solution.
\subsection{A priori estimates}

Let us introduce some notation. For fixed $N$ and given
$v^0,\dots,v^N$ in $H$, we define the {\sl piecewise linear}\/ 
and {\sl backward constant}\/ interpolants of a vector $\{v^i\}_{i=0}^N$
respectively as follows:
\begin{gather}
  \widehat v_\tau(0) := v^0, \quad \widehat v_\tau(t) := a_i(t) v^i +
  \big( 1 - a_i(t) \big) v^{i-1}, \nonumber \\
  \overline v_\tau(0) := v^0, \quad \overline v_\tau(t):= v^i, \qquad
  \text{for} \ \ t \in ((i-1)\tau,i \tau], \ \ i=1, \dots, N ,\label{def_fun}
\end{gather}
where $a_i(t) :=(t-(i-1)\tau)/\tau$ for $t \in ((i-1)\tau,i \tau]$,
$i=1, \dots, N$.
Moreover, let us introduce the translation operator ${\cal T}_\tau$ related to the time
step $\tau$ by setting
\begin{gather}
  ({\cal T}_\tau v)(t) := v(0) \ \ \forall\, t \in [0,\tau) \quad
    \text{and} \ \ ({\cal T}_\tau v)(t)
    := v(t -\tau) \qquad \forall\, t \in [\tau,T]\nonumber \\
  \forall\, v :[0,T] \longrightarrow H. \label{ulisse2}
\end{gather}  
Finally, for $i=1,\dots,N$, we set $\delta
v^i:=(v^i-v\iiu)\tau^{-1}$. It is clear that 
$\overline{\delta v}_\tau=\widehat v_{\tau,t}$ a.e.~in~$[0,T]$.

With this notation \eqref{eqi} takes the form
\beeq{eqtau}
  \alpha\nnu(\widehat u_{\tau,t})
   +B\+\overline u_\tau+\beta\ee(\overline u_\tau)+\epsi\+\overline u_\tau
    =\lambda{\cal T}_\tau \overline u_\tau+\overline f_\tau.
\end{equation}
Let us now prove a technical result, i.e., the extension of 
\eqref{alfaprimo} to $\alpha\nnu$.
\bele\label{coercalfanu}
 There exists $\overline\nu>0$ such that 
 \beeq{alfaprimonu}
   \alpha\nnu'(r)\ge \sigma 
    \qquext{for a.e.~}\,r\in \RR\setminus [2S_-,2S_+],
    ~~\nu\in(0,\overline\nu].
 \end{equation}
\enle
\noindent
Clearly, this also entails the extension of \eqref{alfaprimo2}
to $\alpha\nnu$, with $\sigma/2$ in place of $\sigma$.
\begin{proof}
Let us firstly assume for simplicity $\, \alpha \,$ to be locally Lipschitz continuous. In this case, one clearly has
\begin{equation}\label{label}
(\Id+\nu\alpha)([S_-,S_+])\subset[2S_-,2S_+]\ \ \text{at least for sufficiently small $\nu$.}
\end{equation}
 Then, we have $(\Id+\nu\alpha)'\ge 1+2\nu\sigma$ outside
 $[S_-,S_+]$, whence the resolvent $j\nnu:=(\Id+\nu\alpha)^{-1}$
 satisfies
 \beeq{jnnu} 
   j\nnu'(r)=\frac1{(\Id+\nu\alpha)'(j\nnu(r))}
    \le \frac1{1+2\nu\sigma}
    \qquad\perogni r\not\in[2S_-,2S_+].
 \end{equation}
 Indeed, $r\not\in[2S_-,2S_+]\Rightarrow
 j\nnu(r)\not\in[S_-,S_+]$. From \eqref{jnnu},
 we obtain
 \beeq{annu} 
   \alpha\nnu'(r)=\frac{1-j\nnu'(r)}\nu
   \ge \frac{2\sigma}{1+2\nu\sigma}
    \quad\perogni r\not\in[2S_-,2S_+],
 \end{equation}
 whence the thesis follows for $\nu$ sufficiently small.

By suitably adapting the notations, the latter argument 
still holds for multi-valued graphs $\, \alpha\,$ whenever
the sets $\, \alpha(S_-) \,$ and $\,\alpha(S_+) \,$ are bounded.
Indeed, in this case one can still
establish  $r\not\in[2S_-,2S_+]\Rightarrow
 j\nnu(r)\not\in[S_-,S_+]$ for small $\, \nu$.

On the contrary, when $\, \inf \alpha(S_-) =-\infty\,$ 
(or $\,\sup\alpha(S_+)=+\infty $, or both) relation \eqref{alfaprimonu}
is even easier to check since in this case
\iniziomodificag
$\, \alpha_\nu'(r)= 1/\nu\,$ for all $\, r <2S_-\,$ and
sufficiently small $\nu$.
\finemodificag
Hence \eqref{alfaprimonu} easily follows.
\end{proof}

The same Lemma can be applied to the graph $\beta\ee$. 
More precisely,
one can prove that, for $\eta$ as in \eqref{potpos2b}
and some $\overline\epsi,\overline R>0$,
\beeq{potpos3b}
   W\ee'(r)r\ge\frac{\eta r^2}2
    \quad\perogni r\not \in[-\overline R,\overline R],~
    \epsi\in(0,\overline\epsi].
\end{equation}
Since also the primitive $W\ee$ is defined up to
an integration constant, we can also extend
\eqref{potpos2c}, which takes the form
\beeq{potpos3c}
   W\ee(r)\ge\frac{\eta r^2}4
   \qquad\perogni r \in \RR,~~
     \epsi\le\overline\epsi.
\end{equation}
Let us now prove the a priori estimates. We point out that
all the constants $c$ in this section are 
independent of the approximation parameters
$\epsi$, $\nu$, and $\tau$ as well as of $T$.
Some specific constant will be noted by $c_j$, 
$j\ge 1$.

\paragraph{ First estimate.}~~%
Let us test \eqref{eqi}
by $(u\ii-u\iiu)$
in the scalar product of $H$.
\iniziomodifica
Using the definition of subdifferential
\finemodifica
 and 
the elementary equality
$$
  \tau(v\ii,\delta v\ii)=\frac12
   \big(|v\ii|^2+|v\ii-v\iiu|^2-|v\iiu|^2\big),
$$
we get
\bealo
  & \tau\big(\alpha\nnu(\delta u\ii),\delta u\ii\big)
   +\Phi(u\ii)-\Phi(u\iiu)
   +\betaciapo\ee(u\ii)-\betaciapo\ee(u\iiu)
   +\frac\epsi2|u\ii|^2-\frac\epsi2|u\iiu|^2\\
 \label{conto11}
   & \mbox{}~~~~~
  \le(f\ii,u\ii-u\iiu)
  +\frac\lambda2\big(|u\ii|^2-|u\ii-u\iiu|^2-|u\iiu|^2)
  +c,
\end{align}
where $\betaciapo\ee$ is the primitive of $\beta\ee$
satisfying $\betaciapo\ee(0)=0$. Then (see Lemma
\ref{coercalfanu}),
the first term gives
\beeq{conto11.2}
  \tau\big(\alpha\nnu(\delta u\ii),\delta u\ii\big)
   \ge\frac{\sigma\tau}2|\delta u\ii|^2-c_1\tau,
\end{equation}
with $c_1=0$ if \eqref{newrego1a} holds.
Correspondingly, the 
first term on the \rhs\ of \eqref{conto11}
can be estimated as 
\beeq{conto11.3}
   (f\ii,u\ii-u\iiu)
    \le\frac{\sigma\tau}8|\delta u\ii|^2
    +\frac{2\tau}\sigma|f\ii|^2.
\end{equation}
Hence, let us collect $W\ee(u^i)$ from
its summands in \eqref{conto11},
split it and the $\alpha_\nu$-term 
into two halves, and sum over $i=1,\dots,m$, for
arbitrary $m\le N$. Noting 
also that $\Phi(u^0)<\infty$ (which is a consequence of 
\eqref{hpu0} and of properties of subdifferentials),
using Young's inequality, and taking the maximum
as $m$ ranges in $\{1,\dots,N\}$, we obtain
(for clarity in the $\tau$-notation)
\bealo
  & \frac\sigma8\|\widehat u_{\tau,t}\|_{L^2(0,T;H)}^2
   +\frac12\iTo\alpha_\nu(\widehat u_{\tau,t})\widehat u_{\tau,t}
   +\|\Phi(\overline u_\tau)\|_{L^\infty(0,T)}
   +\frac12\|W\ee(\overline u_\tau)\|_{L^\infty(0,T;L^1(\Omega))}\\
 \label{conto12}
  & \mbox{}~~~~~~~~~~~~~~~
  +\frac\eta{8} \|\overline u\|_{\LIH}^2
  \le\frac2\sigma\|\overline f_\tau\|^2_{L^2(0,T;H)}
   +c_1T+c,
\end{align}
where $c_1$ is the same as in \eqref{conto11.2}
(and is $0$ if \eqref{newrego1a} is fulfilled), while $c$
depends on the initial data and on the integration 
constant chosen such that $W\ee$ fulfills
\eqref{potpos3c} (which also justifies taking
its norm).

We also point out that, if \eqref{newrego2} holds,
then $f\ii=f$ for all $i$. Hence, the first term 
in the \rhs\ of \eqref{conto11} can also be estimated 
in an alternative way. Namely, summing on $i$
we have
\beeq{conto13}
  \sumim(f\ii,u\ii-u\iiu)
   =(f,u^m)-(f,u_0)
   \le\frac\eta{16} |u^m|^2+c,
\end{equation}
where $c$ depends on $\eta,f,u_0$, but neither on $T$ 
nor on $m,N$. In this case, it is 
worth rewriting \eqref{conto12} in a different form,
avoiding use of \eqref{conto11.2} and 
exploiting \eqref{conto13} in place of \eqref{conto11.3}.
Splitting again the $W\ee$-term into two
halves, we then have
\beeq{conto14}
  \iTo\alpha_\nu(\widehat u_{\tau,t})\widehat u_{\tau,t}
   +\|\Phi(\overline u_\tau)\|_{L^\infty(0,T)}
  +\frac12\|W\ee(\overline u_\tau)\|_{L^\infty(0,T;L^1(\Omega))}
  +\frac\eta{16} \|\overline u\|_{\LIH}^2
  \le c,
\end{equation}
where the latter $c$ depends on the constant in \eqref{conto13},
on the integration constant of $W\ee$, and on the initial
data, but is independent of approximation parameters
and of $T$.

\paragraph{Second estimate.}
Let us now test \eqref{eqi} in the scalar
product of $H$ by 
$$
  \tau\delta(B+W'\ee+\epsi)u\ii 
   =(Bu\ii+\beta\ee(u\ii)+\epsi u\ii-\lambda u\ii)
   -(Bu\iiu+\beta\ee(u\iiu)+\epsi u\iiu-\lambda u\iiu).
$$
This gives
\bealo
  & \mbox{}~~~~~~~~~
   \Big(\alpha\nnu\Big(\frac{u\ii-u\iiu}\tau\Big),\tau\delta(B+W'\ee+\epsi)u\ii\Big)\\
  \label{conto20}
  & \mbox{}
   +\big((B+\beta\ee+\epsi)u\ii-\lambda u\iiu,
    \tau\delta(B+W'\ee+\epsi)u\ii\,\big)
   =\big(f\ii,\tau\delta(B+W'\ee+\epsi)u\ii\big).
\end{align}
Then, let us notice that, by monotonicity of 
$\alpha\nnu$ and $\beta\ee$ and $\alpha\nnu(0)=0$,
one has that  
\beeq{conto21}
  \Big(\alpha\nnu\Big(\frac{u\ii-u\iiu}\tau\Big),
    \beta\ee(u\ii)-\beta\ee(u\iiu)+\epsi(u\ii-u\iiu)\Big)
   \ge0.
\end{equation}
Analogously, integrating by parts and using the monotonicity
of $\alpha\nnu$ and $\bbb$ (cf.~\eqref{defib}) 
and the Lipschitz continuity of $\alpha\nnu$,
we get that
\beeq{conto22}
  \Big(\alpha\nnu\Big(\frac{u\ii-u\iiu}\tau\Big),
    B u\ii-B u\iiu\Big)\ge0.
\end{equation}
Next, let us sum over $i=1,\dots,m$, $m\le N$.
Elementary calculations 
permit us to conclude that
\bealo
  & \sumim \big((B+\beta\ee+\epsi)u\ii-\lambda u\iiu,
    \delta(B+W'\ee+\epsi)u\ii\,\big)\\
 \no
   & \mbox{}~~~~~
  =\frac12\big|(B+\beta\ee+\epsi-\lambda)u^m|^2
  -\frac12\big|(B+\beta\ee+\epsi-\lambda)u^0|^2\\
   \label{conto23}
   & \mbox{}~~~~~~~~~~ 
  +\frac12\sumim\big|(B+\beta\ee+\epsi)u\ii-(B+\beta\ee+\epsi)u\iiu\big|^2
  -\frac12\sumim\big|\lambda u\ii-\lambda u\iiu\big|^2.
\end{align}
By collecting the above calculations we have proved that 
\begin{gather}
  \frac12 |(B + W_\epsi' +\epsi)u^m|^2 
\leq \frac12 |(B + W_\epsi' +\epsi)u^0|^2
+ \lambda \sum_{i=1}^m \left(\alpha_\nu\left(\delta u^i\right),
 u^i - u^{i-1}\right) \nonumber\\
+ \frac{\lambda^2}{2}\sum_{i=1}^m|u^i - u^{i-1}|^2
 + \sum_{i=1}^m
(f_i, \tau \delta (B + W_\epsi' +\epsi)u^i)\label{**}
\end{gather}
Furthermore, let us observe that, in case \eqref{newrego2} holds,
i.e.~$f\ii=f$ for all $i$, we get, for $c$ depending on initial
data and on $f$,
\bealo
  \sumim\big(f,\tau\delta(B+W'\ee+\epsi)u\ii\big)
  & =\big(f,(B+W'\ee+\epsi)u^m-(B+W'\ee+\epsi)u^0\big)\\
 \label{conto24} 
  & \le c+\frac14\big|(B+W'\ee+\epsi)u^m\big|^2.
\end{align}  
Finally, by the monotonicity of $B$ and $\beta\ee$, we have 
\begin{gather}
   \frac14\big|(B+\beta\ee+\epsi-\lambda)u^m|^2
   \ge\frac14|Bu^m|^2
   +\frac14|\beta\ee(u^m)|^2
   +\frac{(\epsi-\lambda)^2}4|u^m|^2\nonumber\\
   +\frac{(\epsi-\lambda)}2\big(Bu^m,u^m\big)
   +\frac{(\epsi-\lambda)}2\big(\beta\ee(u^m),u^m\big),\nonumber\\
  \ge\frac18|Bu^m|^2
   +\frac18|\beta\ee(u^m)|^2
   -c_2(1+\lambda^2)|u^m|^2,\label{conto23.2}
\end{gather}
where $c_2$ can be chosen independently of $\epsi$
and $\lambda$.

Collecting \eqref{conto21}--\eqref{conto24}
and taking the maximum for $m\in\{1,\dots,N\}$,
in the $\tau$-notation \eqref{conto20} becomes
\bealo
  & \frac18\|B\overline u_\tau\|^2_{\LIH}
   +\frac18\|\beta\ee(\overline u_\tau)\|^2_{\LIH}\\
 \label{conto25}
  & \mbox{}~~~~~~~
  \le c
   +\lambda\iTo\alpha(\widehat u_{\tau,t})\widehat u_{\tau,t}
   +\frac{\tau\lambda^2}2\big\|\widehat u_{\tau,t}\big\|_{L^2(0,T;H)}^2
   +c_2(1+\lambda^2)\|\overline u_\tau\|_{\LIH}^2.
\end{align}
As before, let us consider also the case where
\eqref{newrego1b} holds in place of \eqref{newrego2}.
Now, in place of \eqref{conto24}, the term with 
$f$ has to be bounded by use
of the following discrete Gronwall-like
Lemma, whose proof can be easily deduced from \cite[Prop. 2.2.1]{Jerome83}.

\bele
 Let $v^i,z^i\in H$, $i=0,1,\dots,N$, be such that
 \beeq{hpgronwdiscr}
  |v^m|^2\le c+c\tau\sum_{i=1}^m \big(z^i,(\delta v)\ii\big)
   \qquad\perogni m=1,\dots,N,
 \end{equation}
for some $\, c>0$.
 Then,  we have
 \beeq{tsgronwdiscr}
   |v^m|\le c'
   \qquad\perogni m=1,\dots,N.
 \end{equation}
for some $\,c'>0\,$ depending only on $\, c,\, v^0,\,z^0$, and $\, \sum_{i=1}^N
\tau |\delta z^i|$.
\enle
%
%
By applying the latter Lemma to \eqref{**} with 
$v^i=(B+W\ee'+\epsi)u^i$ and $z\ii=f\ii$, 
we then arrive again at the estimate \eqref{conto25},
of course with a different value of the 
constant $c$ which now depends in particular on the 
norm of $f$ in $L^1(0,+\infty;H)$.

\paragraph{Comprehensive estimate.}
We consider separately the general case and the specific one
given by 
\eqref{newrego2}. In the general case, we have to multiply
\eqref{conto12} by 
$\max\{4\lambda,8(c_2(1+\lambda^2)+1)\eta^{-1}\}$ and sum the result to
\eqref{conto25}, getting
\bealo
 & \lambda\iTo\alpha_\nu(\widehat u_{\tau,t})\widehat u_{\tau,t}
 +\frac{\lambda\sigma}2\|\widehat u_{\tau,t}\|_{L^2(0,T;H)}^2
   +4\lambda\|\Phi(\overline u_\tau)\|_{L^\infty(0,T)}
   +2\lambda\|W\ee(\overline u_\tau)\|_{L^\infty(0,T;L^1(\Omega))}\\
 \no 
  & \mbox{}~~~~~~~~~~~~~
   +\frac18\|B\overline u_\tau\|^2_{\LIH}
   +\frac18\|\beta\ee(\overline u_\tau)\|^2_{\LIH}
   +\|\overline u_\tau\|_{\LIH}^2\\
 \label{conto31}
  & \mbox{}~~~~~
  \le \frac c\sigma\|\overline f_\tau\|^2_{L^2(0,T;H)}
   +\frac{\tau\lambda^2}2\big\|\widehat u_{\tau,t}\big\|_{L^2(0,T;H)}^2
   +c_3T+c,
\end{align}
where $c_3$, which comes from $c_1$, is 
$0$ if \eqref{newrego1a} holds.

If \eqref{newrego2} holds, we can get something more precise if we
multiply \eqref{conto14} by\linebreak
$\max\{2\lambda,16(c_2(1+\lambda^2)+1)\eta^{-1}\}$ and sum the result to
\eqref{conto25}. This gives
\bealo
 & \lambda\iTo\alpha_\nu(\widehat u_{\tau,t})\widehat u_{\tau,t}
   +2\lambda\|\Phi(\overline u_\tau)\|_{L^\infty(0,T)}
   +\lambda\|W\ee(\overline u_\tau)\|_{L^\infty(0,T;L^1(\Omega))}
  +\|\overline u_\tau\|_{\LIH}^2\\
 \label{conto32} 
  & \mbox{}~~~~~~~~~~~
   +\frac18\|B\overline u_\tau\|^2_{\LIH}
   +\frac18\|\beta\ee(\overline u_\tau)\|^2_{\LIH}
   \le \frac{\tau\lambda^2}2\big\|\widehat u_{\tau,t}\big\|_{L^2(0,T;H)}^2
   +c.
\end{align}
Estimate \eqref{conto31} (or \eqref{conto32}) is the basic
ingredient needed in order to pass to the limit within the approximations.


\subsection{Passage to the limit}
\label{takinglimit}

We remove here first the approximation in $\tau$ and then,
simultaneously, those in $\epsi$ and $\nu$. In order to take
the first limit, let us fix $T>0$. All the convergence properties
listed below are intended to hold up to successive extractions
of subsequences of $\tau\searrow0$, not relabeled.
Let us notice that
the second term in the \rhs\ of \eqref{conto31} 
is estimated, for $\tau$ sufficiently small, by the 
second on the \lhs. 
\iniziomodifica
Moreover, note that a straightforward comparison argument in 
equation \eqref{eqi} entails that $ \alpha\nnu(\widehat u_\tau)$ is
uniformly bounded with respect to $\tau$ in $\LIH$.
\finemodifica
 Hence, for suitable limit functions 
$u,\xi,\eta$ we obtain
\beal{co11}
  \overline u_\tau,\widehat u_\tau\to u
   & \quext{weakly-$*$ in }\,\LIH,\\
 \label{co12}
  \widehat u_{\tau,t}\to u_t
   & \quext{weakly in }\,\LDH,\\
  \label{co13}
 B\overline u_\tau\to \eta
   & \quext{weakly-$*$ in }\,\LIH,\\
  \label{co14}
 \alpha\nnu(\widehat u_\tau)\to \xi
   & \quext{weakly-$*$ in }\,\LIH,
\end{align}
where it is a consequence of standard argument 
the fact that the limits of $\overline u_\tau$ and $\widehat u_\tau$
do coincide.
Moreover, it is clear that \eqref{conto31} entails
\beeq{unifbddE}
  \nor{E(\overline u_\tau)}{L^\infty(0,T)}
   +\nor{E(\widehat u_\tau)}{L^\infty(0,T)}\le c,
\end{equation}
whence, by \eqref{co12}, Lemma~\ref{lemmacompsub},
and \cite[Thm. 1]{Si} it is not difficult to get
%
\beeq{co16}
  \overline u_\tau,\widehat u_\tau\to u
   \quext{strongly in }\,L^\infty(0,T;H).
\end{equation}
More precisely, one has
that
\beeq{co16.3}
  \overline u_\tau(t),\widehat u_\tau(t)\to u(t)
   \quext{strongly in }\,H\quext{for {\sl every~}}\,t\in[0,T]
\end{equation}
(and not just almost everywhere in $(0,T)$).

We shall now pass to the limit in \eqref{eqtau}. By Lipschitz continuity, 
 the term 
$\beta\ee(\ove u_\tau)$ passes to the corresponding limit $\beta\ee(u)$. 
Moreover, it is clear that
${\cal T}_\tau u$ and $\overline f_\tau$
tend, respectively, to $u$ and $f$,
strongly in $\LIH$. Thus, letting $\tau\searrow 0$, 
\eqref{eqtau} gives
\beeq{eqtaulim}
  \xi+\eta+\beta\ee(u)+(\epsi-\lambda)\+u=f,
\end{equation}
and it just remains to identify $\xi=\alpha\nnu(u_t)$ and
$\eta=Bu$.

With this aim, one can, e.g., test again \eqref{eqtau}
by $\widehat u_{\tau,t}$ and integrate 
over $(0, T)$. What one gets is
\bealo
  & \iTT\big(\alpha\nnu(\widehat u_{\tau,t}),\widehat u_{\tau,t}\big)
   =-\Phi(\overline u_\tau(T))+\Phi(\overline u_0)
   -\io\betaciapo\ee(\overline
   u_\tau(T))+\io\betaciapo\ee(\overline u_0)\\
 \label{conto41}
  & \mbox{}~~~~~~~~~~
   +\iTT\big(\overline f_\tau+(\lambda-\epsi)\overline u_\tau,
       \widehat u_{\tau,t}\big).
\end{align}

Consequently,
taking the limsup as (a suitable subsequence of)
$\tau$ goes to $0$, noting that the terms on the \rhs\
can be managed thanks to the strong convergence 
$\overline f_\tau\to f$, \eqref{co12} and 
\eqref{co16}, and using \eqref{co16.3} with the convexity and 
lower semicontinuity of $\Phi$ and $\betaciapo\ee$,
a comparison with the limit equation
\eqref{eqtaulim} permits us to say that
\beeq{conto42}
  \limsup_{\tau\searrow0}
        \iTT\big(\alpha\nnu(\widehat u_{\tau,t}),\widehat u_{\tau,t}\big)
   \le\iTT\big(\xi,u_t\big),
\end{equation}

whence the standard monotonicity argument of, e.g.,
\cite[Prop.~1.1, p.~42]{Ba} yields 
that $\xi=\alpha\nnu(u_t)$ a.e.~in $Q_T$. Finally,
the fact that $\eta= Bu$  is a consequence
of \eqref{co13}, \eqref{co16}, and of the quoted
tool from \cite{Ba}. The passage to the limit in 
$\tau$ is thus completed.

\vspace{2mm}

Let us now briefly detail the passage to the limit with
respect to $\nu$ and $\epsi$. We will send both parameters
to $0$ simultaneously, for brevity. However, 
it is clear that the two limits might also be taken in 
sequence, instead. Thus, let us take a sequence of 
$(\nu_n,\epsi_n)\to (0,0)$ as $n\nearrow +\infty$ and,
 in order to emphasize
dependence on $n$, rename $u_n$ the solution, before 
denoted as $u$, yielded by the $\tau$-limit. It is clear
that, by semicontinuity of norms with respect to weak convergence,
we can take the $\tau$-limit of estimates \eqref{conto31},
\eqref{conto32}, which still hold with respect 
to the same constants $c$ and $c_1$ independent of 
$\epsi$, $\nu$, and $T$. Note that the term with
$\tau\+\lambda^2$ on the \rhs s has now disappeared.
\iniziomodifica
Thus, for suitable limiting functions $u, \eta, \xi$, and $
\gamma$ we have that
\begin{eqnarray}
\iniziomodificag
  \label{cococo} 
 u_n \to u
  & \quext{\iniziomodificag strongly in }\iniziomodificag\,C^0([0,T];H),\\
 %
 \label{co22}
\iniziomodifica
   u_{n,t}\to u_t
   & \quext{weakly in }\,\LDH,\\
  \label{co23}
 Bu_n\to \eta
   & \quext{weakly-$*$ in }\,\LIH,\\
  \label{co24}
 \alpha_{\nu_{n}}(u_n)\to \xi
   & \quext{weakly-$*$ in }\,\LIH,\\
   \label{co25}
   \beta_{\epsi_n}(u_n)\to \gamma
   & \quext{weakly-$*$ in}\, \LIH.
\end{eqnarray}
Now, we pass to the limit as $n\nearrow +\infty$ (hence $\epsi_n, \nu_n \searrow 0$) in
\eqref{eqtaulim} written at level $n$. We obtain 
\begin{equation}\label{eqnlim}
\xi+\eta+\gamma -\lambda u=f,
\end{equation}
and we have to identify $\xi,\eta$, and $ \gamma$ in terms of
$\alpha(u_t), Bu$, and $ \beta(u)$, respectively. 

The identification $\eta \in Bu $ follows immediately from
\eqref{cococo} 
and \eqref{co23}. As for the remaining two inclusions some care 
is needed since the operators themselves are approximations. Indeed, one has that the functionals
$$ L^2(0,T;H) \ni  u \mapsto \int_\Omega \widehat \beta_{\epsi_n} (u(x))\,\dix$$
converges in the sense of Mosco \cite{Att} in $L^2(0,T;H)$ to 
$$L^2(0,T;H) \ni u \mapsto \int_\Omega \widehat \beta (u(x))\,\dix \ \ \text{if} \ \ \beta (u) \in L^1(\Omega) \ \ \text{and} \ \ +\infty \ \ \text{otherwise}. $$
The latter functional convergence, \eqref{cococo}, and \eqref{co25} immediately
give the identification $\gamma\in \beta(u)$  a.e.~in~$Q_T$ via
\cite[Prop.~3.56.c, p.~354 and Prop.~3.59, p.~361]{Att}. 
Moreover, owing to the lower semicontinuity of $\Phi$ and the
convergence \eqref{cococo}, we readily check that
\begin{equation}\label{conto43}
\liminf_{n\nearrow +\infty}\big(\Phi(u_n(T))+\int_{\Omega}\hat{\beta}_{\epsi_n}(u_n(T))\big)\ge 
\Phi(u(T))+\int_{\Omega}\hat{\beta}(u(T)).
\end{equation}
Arguing once again along the lines of \eqref{conto41}, the latter inequality entails in particular that (see \eqref{conto42})
\beeq{conto422}
  \limsup_{n\nearrow +\infty}
        \itt\big(\alpha_{\nu_n}(u_n),u_{n,t}\big)
   \le\itt\big(\xi,u_t\big),
\end{equation}
and the inclusion $\xi\in \alpha(u_t)$ a.e. in $Q_T$ follows from the above-cited results from \cite{Att}.
\finemodifica



The proof of existence is thus concluded. Let us make, anyway,
two final observations. Actually, 
a by-product of
our procedure is that also the limit solution $u$ satisfies
estimates analogous to \eqref{conto31}, \eqref{conto32}
(but {\sl without} the term with $\tau$ on the \rhs). We report,
for completeness, the limit version of \eqref{conto31},
which will be used again in the sequel
\bealo
 & \lambda\iTo\alpha(u_t)\, u_t
  +\frac{\lambda\sigma}2\|u_t\|_{L^2(0,T;H)}^2
   +4\lambda\|\Phi(u)\|_{L^\infty(0,T)}
   +2\lambda\|W(u)\|_{L^\infty(0,T;L^1(\Omega))}\\
 \no
  & \mbox{}~~~~~~~~~~~~~~~~~~~~~~~~
   +\frac18\|B\overline u\|^2_{\LIH}
   +\frac18\|\beta(u)\|^2_{\LIH}
   +\|u\|_{\LIH}^2\\
  \label{conto31lim}
 & \mbox{}~~~~~~~~
   \le c\big(1+\|f\|^2_{L^2(0,T;H)}\big)+c_3T,
\end{align}
where $c_3=0$ if 
\eqref{newrego1a} holds and still $c,c_3$ are
independent of $T$.\\
\iniziomodifica
The second observation is that, thanks to the convergence \eqref{co23}, we can prove
that 
\begin{equation}\label{co26}
\Phi(u_n(t))\to \Phi(u(t)) \quad\text{for any }\, t\in [0,T].
\end{equation}
Actually, the definition of subdifferential written for 
$u_n$ gives that
\begin{equation}\label{conto44}
\Phi(u_n(t))\le (Bu_n(t),u_n(t)-u(t))+\Phi(u(t)) \quad\text{for {\sl any} }\, t\in [0,T]
\end{equation}
\iniziomodificag
(and not only almost everywhere, see the Proof of
Proposition~\ref{teosep} just below).
\finemodificag
Thus, by taking the limsup in \eqref{conto44} and recalling \eqref{cococo}, \eqref{co23},  and the lower 
semicontinuity of $\Phi$, we readily get \eqref{co26}. In the specific case in which 
$\Phi(v)=\frac{1}{p}\|\nabla v\|_{p}^{p}$ (hence $Bv$ is the $p$-laplacian), the convergence
\eqref{co26} entails the convergence of $u_n$ in $W^{1,p}(\Omega)$ (recall that $p>1$ and thus 
 $W^{1,p}(\Omega)$ is uniformly convex).
\finemodifica


\section{Separation property and uniqueness}
\label{secseparuniq}

Henceforth, let us
denote by {\sl solution} any function $u$ satisfying
\eqref{eqn}--\eqref{iniz} in the sense and with the 
regularity made precise in Theorem~\ref{teoexist}.

\paragraph{ Proof of Proposition~\ref{teosep}.}
Let us first prove that, for any solution $u$,
one has $u(t)\in \dom_H B$ for all, and not just 
a.e., $t\in[0,T]$. Indeed, by \eqref{regu1}, 
$u$ lies in $C^0([0,T];H)$. Thus, assuming by contradiction
that, for some $t$, $u(t)\not\in\dom_H B$, we can approximate
$t$ by a sequence $\{t_n\}$ such that $u(t_n)\in\dom_H B$
for all $n$. Since \eqref{reguunif} holds, we can 
assume $\{Bu(t_{n})\}$ bounded in $H$. Hence, 
extracting a subsequence $\{n_k\}$
such that $u(t_{n_k})\to u(t)$ strongly in $H$ and 
$Bu(t_{n_k})\to \widetilde B$ weakly in $H$,
by maximal monotonicity of $B$ we have that
${\widetilde B}=Bu(t)$, whence the assert 
follows.

The same argument can be applied in order to check that
$\beta^0(u(t))\in H$ for {\sl all}\/ $t\in[0,T]$.
As a consequence, by semicontinuity
of norms with respect to weak convergence, we have 
more precisely that
\beeq{nonsoloqo}
  |u(t)|+|Bu(t)|+|\beta^0(u(t))|\le c
   \quad\perogni t\in[0,T],
\end{equation}
where $c$ is the same constant as in \eqref{reguunif}
and, in particular, does not depend on $T$ if 
either \eqref{newrego1a}--\eqref{newrego1b} or
\eqref{newrego2} hold.
Consequently, 
thanks to \eqref{eqsep}-\eqref{eqsepbis},
there exists one constant $\delta>0$,
depending only on $c$ in  
\eqref{nonsoloqo} through the
function $\gamma$, such that 
\beeq{hoelder}
  |u(x_1,t)-u(x_2,t)|\le \delta|x_1-x_2|^\nu
   \quad\perogni x_1,x_2\in\barO,~t\in[0,T],
\end{equation}

We can now prove the right inequality in \eqref{sepsu} 
in the case when $I$ is right-bounded (if $I$ is not
right-bounded, the inequality is trivial). The proof
of the left inequality is analogous, of course. Set
\beeq{defirho}
   \rho:=\delta^{-1/\nu}|r_1-\barr|^{1/\nu}
\end{equation}
and assume, by contradiction, that there exist
$\bart>0$, $\barx\in\barO$ such that 
$u(\barx,\bart)=\barr$. By \eqref{hoelder},
it is clear that
\beeq{hoelder2}
  |u(x,\bart)-\barr|
   =|u(x,\bart)-u(\barx,\bart)|
   \le \delta |x-\barx|^\nu
   \le |r_1-\barr|
  \quad \perogni x\in\barO\cap B(\barx,\rho).
\end{equation}

Since the value of $u$ cannot exceed $\barr$, \eqref{hoelder2}
entails that $u(x,\bart)\ge r_1$ for all 
$x\in\barO\cap B(\barx,\rho)$. Then, by 
\eqref{crescibeta},
\beeq{hoelder3}
  \beta^0(u(x,\bart))
   \ge \frac c{\big(\barr-u(x,\bart)\big)^\kappa}
   =\frac c{\big(u(\barx,\bart)-u(x,\bart)\big)^\kappa}
  \quad \perogni x\in\barO\cap B(\barx,\rho).
\end{equation}
By taking squares, integrating in space, and using \eqref{nonsoloqo}
and \eqref{hoelder2}, we obtain
\bealo
  c
   & \ge \io (\beta^0)^2(u(x,\bart))\,\dix
     \ge  \int_{\Omega\cap B(\barx,\rho)}  (\beta^0)^2(u(x,\bart))\,\dix\\
 \label{contradic}
   & \ge \int_{\Omega\cap B(\barx,\rho)} 
      \frac c{|u(x,\bart)-u(\barx,\bart)|^{2\kappa}}\,\dix
    \ge \int_{\Omega\cap B(\barx,\rho)}
      \frac{c\delta^{-2\kappa}}{|x-\barx|^{2\kappa\nu}}\,\dix.
\end{align}
Now, since $\Omega$ is a smooth set 
(here Lipschitz would be enough), there exists $c_\Omega>0$ such
that, for any sufficiently small $r>0$, $\Omega\cap B(\barx,r)$
measures at least $c_\Omega r^d>0$. Recalling \eqref{kappa}, 
this entails that the latter integral in \eqref{contradic} is $+\infty$,
yielding a contradiction.

This means that no solution $u$ can ever attain the value $\barr$. 
However, in order to show \eqref{sepsu}, we have to be 
more precise. We actually claim that, if
${\cal F}\subset H$ is any set such that
\beeq{boundcalF}
  \esiste\delta'>0:~~\calF\subset\calG(\delta'), 
   \quext{where }~
    \calG(\delta'):=\big\{v\in H:|v|+|Bv|+|\beta^0(v)|\le \delta'\big\},
\end{equation}
then there exists $\barr^*$ such that $v(x)\le\barr^*$ for all $v\in
{\cal F}$ and $x\in\barO$. Applying this to the family
$\calF=\{u(t)\}_{t\in[0,T]}$, we clearly get the upper inequality
in \eqref{sepsu}, as desired. In addition,
it is a by-product of the argument that $\barr^*$ is
independent of $T$ if such is $c$ in \eqref{reguunif}, i.e., 
if either \eqref{newrego1a}--\eqref{newrego1b} or
\eqref{newrego2} hold. 

To prove the claim, let us proceed by contradiction. Namely,
suppose that there are $\{v_n\}\subset \calF$, $\{x_n\}\subset\barO$
such that $v_n(x_n)\nearrow\barr$. Since $\calG(\delta')$ is bounded in 
$C^{0,\nu}(\barO)$ by \eqref{eqsep}-\eqref{eqsepbis}, then
we can extract a subsequence ${n_k}$ such that $x_{n_k}\to \barx\in\barO$ and
$v_{n_k}\to v$ uniformly. Thus, $v(\barx)=\barr$. But this 
is impossible, since it can be easily seen that 
that also $v\in\calG(\delta')$;
hence, for the first part of the argument
$v$ can never attain the value $\barr$.

\noindent%
{\bf Proof of Theorem~\ref{uniq}.}~~%
Assume by contradiction there exist two 
solutions $u_i$, $i=1,2$. 
Then, writing \eqref{eqn} for $i=1,2$, taking the difference,
and testing it by $(u_1-u_2)_t$, we get, for a.e.~$t>0$,
\beeq{conto50}
  \io \big(\alpha((u_1)_t)-\alpha((u_2)_t)
   +B(u_1-u_2)
   +W'(u_1)-W'(u_2)\big)(u_1-u_2)_t=0.
\end{equation}
Next, we note that
\beeq{conto51}
  \io B(u_1-u_2)(u_1-u_2)_t
   =\frac12\ddt\io B(u_1-u_2)(u_1-u_2)
\end{equation}
and, by ({\bf H2}) and \eqref{newrego1a},
\beeq{conto52}
  \io \big(\alpha((u_1)_t)-\alpha((u_2)_t)\big)(u_1-u_2)_t
   \ge 2\sigma\io \big|(u_1)_t-(u_2)_t\big|^2.
\end{equation}
Finally, by \eqref{sepsu} and 
\eqref{betalip}, it is clear that
\beeq{conto53}
  \io \big(W'(u_1)-W'(u_2)\big)(u_1-u_2)_t
   \le \frac\sigma2 \big|(u_1)_t-(u_2)_t\big|^2
   +c|u_1-u_2|^2,
\end{equation}
with the last $c$ depending only on $\, \sigma$, on the constant $c$ in 
\eqref{reguunif}, and on the Lipschitz constant
of $W'$ in the interval $[\subr_*,\barr^*]$.

Taking the integral of \eqref{conto50} on $(0,t)$ 
and exploiting the relation
\beeq{conto54}
  |(u_1-u_2)(t)|^2
   \le  t \big\|(u_1)_t-(u_2)_t\big\|_{L^2(0,t;H)}^2
\end{equation}
the conclusion follows immediately by taking
\eqref{conto51}--\eqref{conto53} into account
and applying Gronwall's Lemma.


\section{Long-time behavior}
\label{seclongtime}

{\bf Proof of Theorem~\ref{longtime}}.~~%
Let $\{t_n\}$ be fixed in such a way that $\, t_n \nearrow +\infty$. Then, property \eqref{coinfty1} for some suitable (not relabeled) subsequence
is a direct consequence of bound \eqref{reguunif} and 
the precompactness in $H$ of the trajectory from Lemma~\ref{lemmacompsub}.

In order to conclude, we have to show that the limit $\ui$ solves
the stationary problem \eqref{prostaz}. To see this, let us 
first note that, if \eqref{newrego1b} holds, then 
$f(t)$ tends to $0$ strongly in $H$ as $t\to\infty$. Then,
we can consider the sequence of Cauchy problems for $\, f_n(\cdot):=f(\cdot+t_n)$
\beeq{Cauchyn}
  \begin{cases}
    \alpha((u_n)_t)+Bu_n+W'(u_n)= f_n & \quext{in }\,H\\
    u_n(0)=u(t_n),
  \end{cases}
\end{equation}
and it is clear that $u_n(t):=u(t+ t_n)$ solves \eqref{Cauchyn},
e.g., for $t\in(0,1)$. Moreover, by \eqref{reguunif}, where
$c$ is now independent of $t$, we have
\begin{eqnarray}
  u_n&\to& \util\quext{strongly in }\/C^0([0,1];H),\label{colim1}\\
Bu_n &\to& \widetilde B\quext{weakly-$\ast$ in}\ \ L^\infty(0,1;H), \label{colim2}\\
W'(u_n)&\to& \widetilde w\quext{weakly-$\ast$ in}\ \ L^\infty(0,1;H),\label{colim22}
\end{eqnarray}
for suitable limit functions $\util,\widetilde B,\widetilde w$.
Due to the standard monotonicity argument \cite[Lemma~1.3, p.~42]{Ba},
this immediately yields 
\beeq{identlim}
\widetilde B = B \tilde u, \quad \widetilde w = \beta(\util)-\lambda \util
   \quext{a.e.~in }\/\Omega\times(0,1).
\end{equation}
Furthermore, by \eqref{conto31lim}, which now holds with 
$c_3=0$, 
\beeq{colim3}
  (u_n)_t\to 0\quext{strongly in }\/L^2(0,1;H),
\end{equation}
whence $\util$ is constant in time. Since $\util(0)=\ui$ by
\eqref{coinfty1} and the Cauchy condition in \eqref{Cauchyn},
we readily conclude that $\util\equiv\ui$ for 
a.e.~$t\in(0,1)$. Finally, still from \eqref{reguunif} 
we have 
\beeq{colim4}
  \alpha((u_n)_t)\to \widetilde\alpha
   \quext{weakly-$\ast$ in }\/L^\infty(0,1;H),
\end{equation}
where actually $\widetilde\alpha\equiv 0$ a.e.~in $(0,1)$ thanks to 
\eqref{colim3}, \cite[Lemma~1.3, p.~42]{Ba}, and \eqref{alpha00}.
This completes the proof of relation \eqref{prostaz} and of the
Theorem.

\vspace{2mm}

\noindent%
{\bf Proof of Theorem~\ref{Loj}}.~~%
Let us consider the 
non-autonomous case when \eqref{newrego1b} holds, the
situation where we have \eqref{newrego2} being simpler. 
We argue 
along the lines of \cite[Sec.~3]{CJ}. However, 
due to the presence of the nonlinearity $\alpha$, 
our proof presents further technical 
complications. Let $\ui$ be an element of 
the $\omega$-limit and note first that, 
by precompactness
(cf.~Remark~\ref{rempreco}), $\ui$ is the limit
in $C(\barO)$ of some sequence $\{u(t_n)\}$.
Thus, by \eqref{uisepar}, at least for $n$ sufficiently
large, we have that 
\beeq{unsepar}
    \igiu<u_n(x)<\isu \quad\perogni x\in\Omega.
\end{equation}
This justifies the application of the \L ojasiewicz-Simon inequality,
since $u_n$ eventually ranges in the interval where $W'$ is
analytic. In particular, the barriers at the extrema of 
$W$ are excluded even in case \eqref{sepsu} does not hold.

Then, similarly with \cite{CJ}, we can set (but note that we use
here the norm in $H$ instead of that in $V$ \iniziomodifica(cf.~the regularity 
of $u_t$ and  \eqref{A1} below)\finemodifica,
in agreement with version \eqref{Lojfeireisl}
of the \L ojasiewicz-Simon inequality):
\beeq{jen1}
  \Sigma:=\big\{
    t>0:\|u(t)-\ui\|_V\le \epsilon/3\big\}.
\end{equation}
Clearly, $\Sigma$ is unbounded. Next, 
for $t\in\Sigma$, we put
\beeq{jen2}
  \tau(t):=\sup\big\{
    t'\ge t:\sup_{s\in[t,t']}\|u(s)-\ui\|_V
     \le \epsilon\big\},
\end{equation}
where, by continuity, $\tau(t)>t$ for all
$t\in\Sigma$. Let us fix  $t_0\in\Sigma$ and
divide ${\cal J}:=[t_0,\tau(t_0)]$ into two subsets:
\beal{A1}
  & A_1:=\Big\{
   t\in\calJ:|u_t(t)|\ge\Big(\int_t^{\tau(t_0)}
   |f(s)|^2\,\dis\Big)^{1-\theta}\Big\},\\
 \label{A2}
  & A_2:=\calJ\setminus A_1.
\end{align}
Letting now
\beeq{defiE0}
  \Phi_0(t):=E(u(t))-E(\ui)
   +\frac1\sigma\int_t^{\tau(t_0)}|f(s)|^2\,\dis,
\end{equation}
exploiting assumption \eqref{alfapotenza}
and H\"older's inequality and making a comparison
 in \eqref{eqn},
it is not difficult to see that
\beeq{jen3} 
  \Phi_0'(t)\le
   -\frac\sigma2|u_t(t)|^2
   -\kin\|u_t(t)\|^{\pin+1}_{L^{\pin+1}(\Omega)}
   -\frac1{2\sigma}|f(t)|^2.
\end{equation}
Note that $\Phi_0$ is absolutely continuous
 thanks to \cite[Lemme~3.3, p.~73]{Br}.
Thus, we have \cite[(3.2)]{Je}
\beeq{jen4} 
  \ddt\big(|\Phi_0|^{\theta}\sign\Phi_0\big)(t)\le
   -\theta|\Phi_0(t)|^{\theta-1}
   \Big(
   \frac\sigma2|u_t|^2
   +\kin\|u_t\|^{\pin+1}_{L^{\pin+1}(\Omega)}
   +\frac1{2\sigma}|f|^2\Big)(t).
\end{equation}
Noting that \eqref{Lojfeireisl} 
can be applied and making a further comparison of terms
in \eqref{eqn}, we have for any such $t_0$ and $t\in A_1$
\bealo
  |\Phi_0(t)|^{1-\theta}
   & \le|E(u(t))-E(\ui)|^{1-\theta}
    +\Big|\frac1\sigma\int_t^{\tau(t_0)}|f(s)|^2\,\dis\Big|^{1-\theta}\\
 \no
  & \le c_\ell\|\alpha(u_t(t))\|_{V^*}
  +c_\ell\|f(t)\|_{V^*}
    +\Big|\frac1\sigma\int_t^{\tau(t_0)}|f(s)|^2\,\dis\Big|^{1-\theta}\\
  & \le c\Big(|u_t(t)|
   +\|u_t(t)\|^{\qin}_{L^{{{\TeXchi\qin}}}(\Omega)}
   +|f(t)|\Big),\label{jen5}
\end{align}
where we used the continuous embeddings
$H\subset L^{\TeXchi}(\Omega)\subset V^*$ and the last constant $c$ 
also depends on $\sigma'$ and $\lin$ in \eqref{alfapotenza}.

Thus, being $\chi\qin\le\pin+1$ by \eqref{relapq}
and $\qin\ge(\pin+1)/2$ since $\qin\ge\pin\ge1$,
from \eqref{reguunif} and \eqref{jen5} we have that
\beeq{jen6}
  |\Phi_0(t)|^{\theta-1}
   \ge c\Big(|u_t(t)|
   +\|u_t(t)\|^{{(\pin+1)}/2}_{L^{\pin+1}(\Omega)}
   +|f(t)|\Big)^{-1}.
\end{equation}
Collecting now \eqref{jen3}, \eqref{jen5}, 
and \eqref{jen6}, \eqref{jen4}
gives 
\beeq{jen7} 
   \Big(|u_t(t)|
   +\|u_t(t)\|^{{(\pin+1)}/2}_{L^{\pin+1}(\Omega)}
   +|f(t)|\Big)
   \le -c
  \ddt\big(|\Phi_0|^{\theta}\sign\Phi_0\big)(t),
\end{equation}
whence, integrating over $A_1$ and exploiting that
$\Phi_0$ is a decreasing function (cf.~\eqref{jen3}), we 
get that $|u_t|$ is integrable over $A_1$.

>From this point on, the argument proceeds exactly
as in \cite{CJ}. Namely, by definition of $A_2$ and 
\eqref{decadef} and possibly taking some smaller $\,\theta$, 
one immediately 
gets that $|u_t|$ is integrable over $A_2$ and 
hence on $\calJ$. This permits to show by a simple 
contradiction argument that $\tau(t_0)=\infty$ as
$t_0\in\Sigma$ is sufficiently large. This entails
$u_t\in L^1(t_0,+\infty;H)$, whence the convergence of
the whole trajectory to $\ui$ follows, as desired.

Finally, let us just briefly outline the changes
to be done when, instead, \eqref{newrego2}
holds.
The most significant difference is that now it is convenient to
include $f$ into the energy, setting, for $v\in H$,
$E_f(v):=E(v)-(f,v)_H$. Then, it is clear that the 
\L ojasiewicz-Simon inequality \eqref{Lojfeireisl}
still holds in the form
\beeq{Lojfeireislf}
  |E_f(v)-E_f(\ui)|^{1-\theta}
   \le c_\ell\|Bv+W'(v)-f\|_{V^*}.
\end{equation}
At this point, the proof is performed 
similarly, provided that one defines $\Phi_0$ in~\eqref{defiE0}
with $E_f$ in place of $E$ and without the integral term.
Moreover, one directly gets the integrability
of $|u_t(t)|$ on $\calJ$ and no longer needs to split
it into the two subsets $A_1$ and $A_2$. The details
of the argument are left to the reader.
\beos 
 Let us note that, if \eqref{newrego1b} and \eqref{decadef}
 hold, then it is also possible to estimate
 the decay rate of solutions
 as in \cite{GPS1,GPS2}.
 Namely, one can prove (cf., e.g., \cite[(3.7)]{GPS2}) that
 \beeq{coinfty3}
   |u(t)-\ui|\le c\+t^{-\mu} \quad\perogni t>0,
 \end{equation}
 where $\mu>0$ depends on $\theta,\xi$
 and $c$ depends only on data (and in particular not on time).
\eddos









\end{document}